\documentclass[12pt]{amsart}
\usepackage{amssymb}
\usepackage{graphicx}
\newtheorem{Thm}{Theorem}[section]
\newtheorem{Lem}[Thm]{Lemma}

\newtheorem{Prop}[Thm]{Proposition}
\newtheorem{Rem}[Thm]{Remark}
\newtheorem{Cor}[Thm]{Corollary}

\def\bbarp{\overline{\overline{P}}}
\def\bbarl{\overline{\overline{\LL^+}}}
\def\rank{{\rm{rank}}}
\def\cl{{\rm{cl}}}
\def\Aa{{\mathcal{A}}}
\def\BB{{\mathcal{B}}}

\def\DD{{\mathcal{D}}}
\def\FF{{\mathcal{F}}}
\def\LL{{\mathcal{L}}}
\def\RR{{\mathcal{R}}}

\def\vap{\varphi}

\def\De{\Delta}
\def\Om{\Omega}
\def\C{{\mathbb{C}}}
\def\R{{\mathbb R}}
\def\Z{{\mathbb{Z}}}
\def\Q{{\mathbb{Q}}}

\def\wh{\widehat}
\def\no={\,{\,|\!\!\!\!\!=\,\,}}
\def\sbst{\subseteq}
\def\oli{\overline}

\def\({\left(}
\def\){\right)}
\def\wt{\widetilde}

\def\ssm{\smallsetminus}
\def\spst{\supseteq}
\def\beq{\begin{eqnarray}}
\def\eeq{\end{eqnarray}}
\def\B{b}

\begin{document}

\title[Geometrically constructed homology bases] {Geometrically constructed  bases for
homology of partition lattices of \\  types $A$, $B$ and $D$}

\author[Bj\"orner]{Anders Bj\"orner$^1$}
\address{Royal Institute of Technology, Department of Mathematics,\newline
S-100~44~Stockholm, Sweden}
\email{bjorner@math.kth.se}
\author[Wachs]{Michelle Wachs$^2$}
\address{University of Miami, Department of Mathematics,\newline Coral Gables, FL
331~24, USA}
\email{wachs@math.miami.edu}
\date{\today}

\thanks{1. Supported in part by G\"oran Gustafsson Foundation for Research in Natural Sciences
and Medicine.}
\thanks{2. Supported in part by National Science Foundation grants DMS-9701407 and
DMS-0073760.}

\maketitle
\begin{center} {\it Dedicated to Richard Stanley on the occasion of his 60th birthday} \end{center}

\begin{abstract}{We use the theory of hyperplane arrangements to construct 
natural bases  for the homology of partition lattices of types $A$, $B$ and $D$.  This
extends and explains the ``splitting basis'' for the homology of the partition lattice
given in \cite{Wa96}, thus answering a question asked by R.
Stanley.  

More explicitly,  the following general technique is presented and utilized.  Let
$\Aa$ be a central and essential hyperplane arrangement in
$\R^d$. Let $R_1,\dots,R_k$ be the bounded regions of a generic hyperplane section of
$\Aa$. We show that there are induced polytopal cycles $\rho_{R_i}$ in the homology of
the proper part $\oli{L}_\Aa$ of the intersection lattice such that
$\{\rho_{R_i}\}_{i=1,\dots,k}$ is a basis for $\wt{H}_{d-2} (\oli{L}_\Aa)$.
This geometric method for constructing combinatorial homology bases is applied to the Coxeter
arrangements of types $A$, $B$ and $D$, and to some interpolating arrangements.}
\end{abstract}

\section{Introduction}

In \cite{Wa96} Wachs constructs a basis for the homology of the partition lattice 
$\Pi_n$ via
a certain natural ``splitting'' procedure for permutations. This basis has
very favorable properties with respect to the representation of the symmetric group $S_n$ on
$\wt{H}_{n-3}(\Pi_n, \C)$, a representation that had earlier been
studied by Stanley \cite{St82}, Hanlon \cite{H81} and many others.  It 
also is the shelling basis for a certain  EL-shelling of the partition
lattice given in \cite[Section 6]{Wa96}.  This basis has connections to the free Lie
algebra as well; see \cite{Wa98}.  

We now give a brief description of the splitting basis of \cite{Wa96}.  For each $\omega \in S_n$,
let
$\Pi_\omega$ be the subposet of $\Pi_n$ consisting of partitions obtained by splitting $\omega$. 
In Figure 1 the subposet $\Pi_{3124}$ of $\Pi_4$ is shown.  Each poset $\Pi_\omega$ is
isomorphic to the face lattice of an $(n-2)$-dimensional simplex.  Therefore $\Delta(\oli
\Pi_\omega)$, the order complex of the proper part of $\Pi_\omega$, is an $(n-3)$-sphere embedded
in 
$\Delta(\oli
\Pi_n)$, and hence it determines a fundamental cycle $\rho_\omega \in \tilde H_{n-3}(\oli \Pi_n)$.
In \cite{Wa96} it is shown that a certain subset of  $\{\rho_\omega| \omega \in S_n\}$ forms a
basis for
$\tilde H_{n-3}(\oli
\Pi_n)$;  namely, the set of all $\rho_\omega$ such that $\omega$ fixes $n$.

\vspace{.2in}
\begin{center}
\includegraphics[width=9cm]{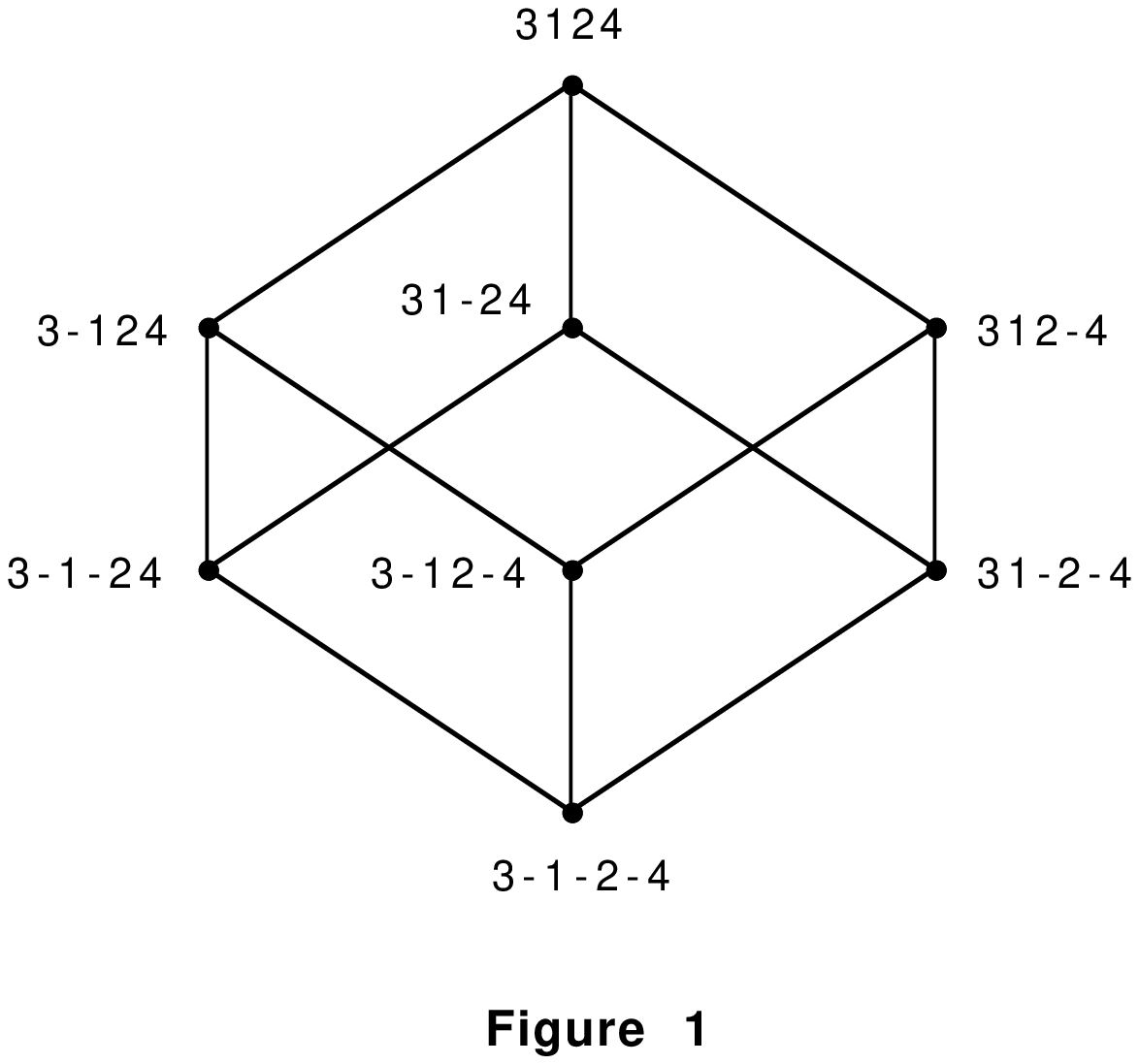}
\end{center}

The  partition lattice is the intersection lattice of
the type $A$ Coxeter arrangement.  
  The original
motivation for this paper was to explain and generalize to other Coxeter groups, the
splitting basis for $\Pi_n$. Taking a geometric point of view we give
such an explanation, which then leads to the construction of
``splitting bases'' also for the intersection lattices of Coxeter arrangements
of types $B$ and $D$ and of  some interpolating arrangements.  Our technique is general in that it
gives a way to construct a basis for the homology of the intersection lattice of any real
hyperplane arrangement.

The intersection lattice of the type $B$ Coxeter arrangement is isomorphic to the  signed
partition lattice   $\Pi^B_n$.  Its elements are signed 
partitions of
$\{0,1,\dots,n\}$; that is  partitions of
$\{0,1,\dots,n\}$ in which any element but the smallest element of each nonzero block can be
barred.  For each element
$\omega$ of the hyperoctahedral group
$B_n$, we  form a subposet $\Pi_\omega$ of $\Pi^B_n$ consisting of all signed partitions obtained
by splitting the signed permutation $\omega$.  In Figure 2 the subposet $\Pi_{\bar 2 3 1}$ of
$\Pi^B_3$ is shown. Just as for type
$A$, it is clear that each subposet $\Pi_\omega$ 
 determines a fundamental cycle $\rho_\omega$ in
$\tilde H_{n-2}(\oli {\Pi^B_n})$. It is not clear, however, that the elements $\rho_\omega$,
$\omega
\in B_n$, generate $\tilde H_{n-2}(\oli{\Pi^B_n})$; nor is it clear how one would select
cycles
$\rho_\omega$ that form a basis for $\tilde H_{n-2}(\oli{\Pi^B_n})$.  Our geometric
technique enables us to identify a basis whose elements are those $\rho_\omega$ for which  the
right-to-left maxima of
$\omega$ are unbarred.

\vspace{.2in}
\begin{center}
\includegraphics[width=9cm]{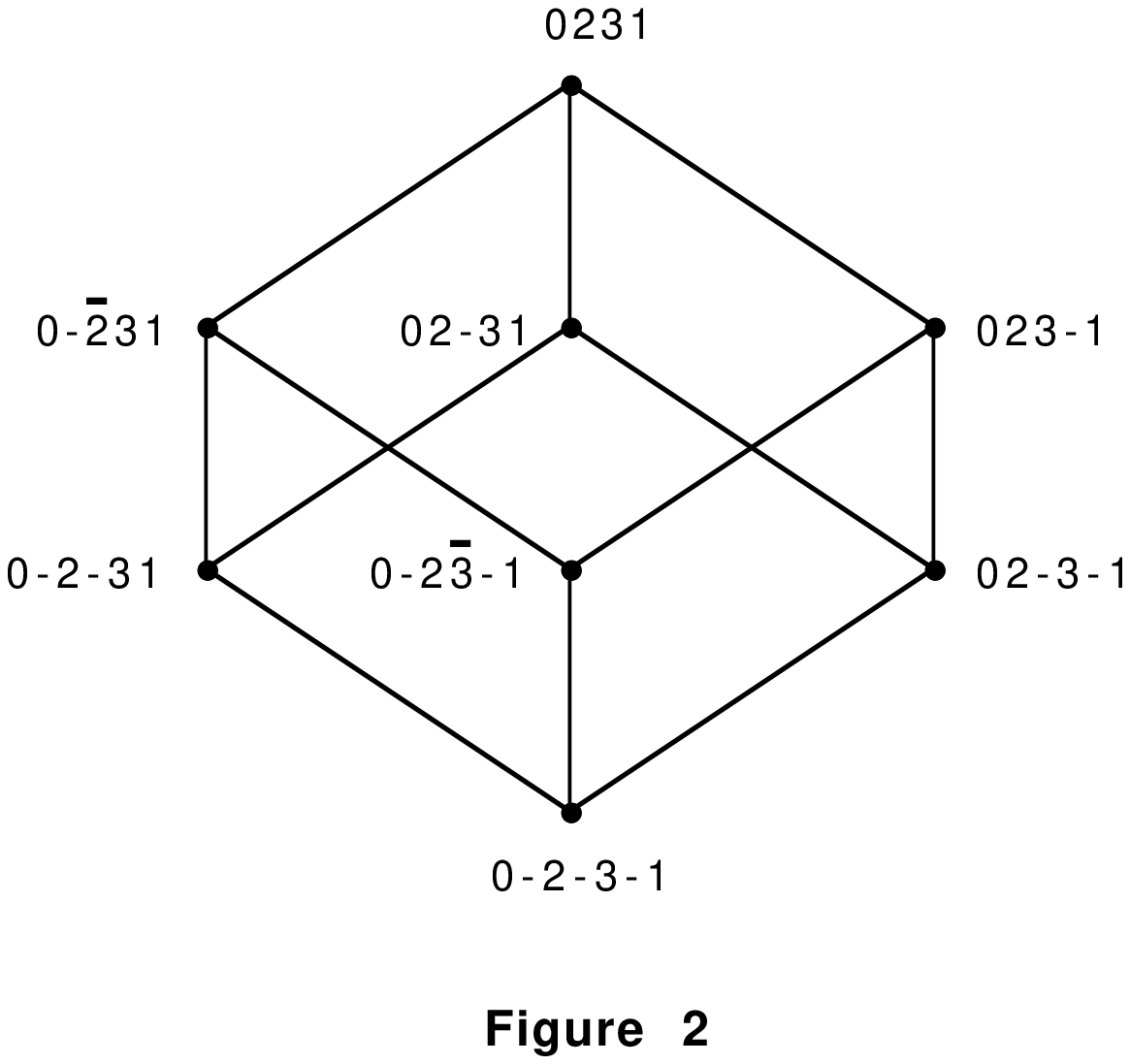}
\end{center}

We will now give a somewhat more detailed description of the content
of the paper. The proper setting for our discussion is that of real
hyperplane arrangements, or (even more generally) oriented matroids.

Let $\Aa$ be an arrangement of linear hyperplanes in $\R^d$. We assume
that $\Aa$ is essential, meaning that $\bigcap \Aa :=
\bigcap_{H\in\Aa} H=\{0\}$. The
intersection lattice $L_\Aa$ is the family of intersections of
subarrangements $\Aa'\sbst\Aa$, ordered by reverse inclusion. It is a
geometric lattice, so it is known from a theorem of Folkman \cite{Fo66}
that $\wt{H}_{d-2} (\oli{L}_\Aa)\cong \Z^{|\mu_L (\hat 0, \hat 1)|}$
and $\wt{H}_i (\oli{L}_\Aa)=0$ for all $i\neq d-2$, where
$\oli{L}_\Aa=L_\Aa-\{\hat 0, \hat 1\}$. In fact, the order complex
$\De(\oli{L}_\Aa)$ has the homotopy type of a wedge of
$(d-2)$-spheres.

There are many copies of the Boolean
lattice $2^{[d]}$ (or equivalently, the face lattice of the
$(d-1)$-simplex) embedded in every geometric lattice of length $d$. Each such Boolean subposet
determines a fundamental cycle in homology.  In
\cite{Bj82} Bj\"orner gives a combinatorial method  for constructing homology bases using
such Boolean cycles. This method, which in its simplest version is based on the so
called ``broken circuit'' construction from matroid theory, is applicable to {\em
all\,} geometric lattices (not only to intersection lattices of
hyperplane arrangements).  Although the cycles in the splitting basis are Boolean,
the basis does not arise from the broken circuit construction. It
turns out that the splitting basis does arise from the geometric
construction in this paper. 

There is a natural way to associate polytopal cycles in the intersection lattice $L_\mathcal A$
with regions of the arrangement $\mathcal A$.  These cycles are not necessarily Boolean.  They are
fundamental cycles determined by face lattices  of convex $(d-1)$-polytopes embedded in
$L_\mathcal A$.
 We show that these cycles generate the 
homology of $\oli L_\Aa$.  Moreover, we present a  way of identifying those regions whose
corresponding cycles form a basis.  Here is a short and non-technical statement of the method.

Let $H$ be an affine hyperplane in $\R^d$ which is generic with
respect to $\Aa$. The induced affine arrangement $\Aa_H=\{H\cap K\mid
K\in\Aa\}$ in $H\cong \R^{d-1}$ will have certain regions that are
bounded. Each bounded region $R$ is a convex $(d-1)$-polytope in $H$
and it is easy to see that a copy of its face lattice sits embedded in
$L_\Aa$. Briefly, every face $F$ of $R$ is the intersection of the
maximal faces containing it, and so $F$ can be mapped to the
intersection of the linear spans (in $\R^d$) of these maximal faces,
which is an element of $L_\Aa$. Thus, we have a cycle $\rho_R \in
\wt{H}_{d-2}(\oli{L}_\Aa)$ for each bounded region $R$. A main
result (Theorem~\ref{T2}) is that these cycles $\rho_R$, indexed by the
bounded regions of $\Aa_H$, form a basis for $\wt{H}_{d-2}(\oli{L}_\Aa)$.

The regions of a Coxeter arrangement are simplicial cones that correspond bijectively to the
elements  of the 
Coxeter group.  When the geometric method is applied to the intersection
lattice of any Coxeter arrangement, the cycles in the resulting basis are Boolean and
are indexed by the elements of the Coxeter group that correspond to the bounded regions of a
generic affine slice.  For type A,  when the generic affine hyperplane $H$ is chosen appropriately
one gets the splitting  basis consisting of cycles
$\rho_\omega$  indexed by the permutations
$\omega$ that fix $n$.  In Figure 3  the intersection of the  Coxeter arrangement
$\mathcal A_3$ with $H$  is shown.  The bounded regions are labeled by their
corresponding permutation.

\vspace{.2in}
\begin{center}
\includegraphics[width=9cm]{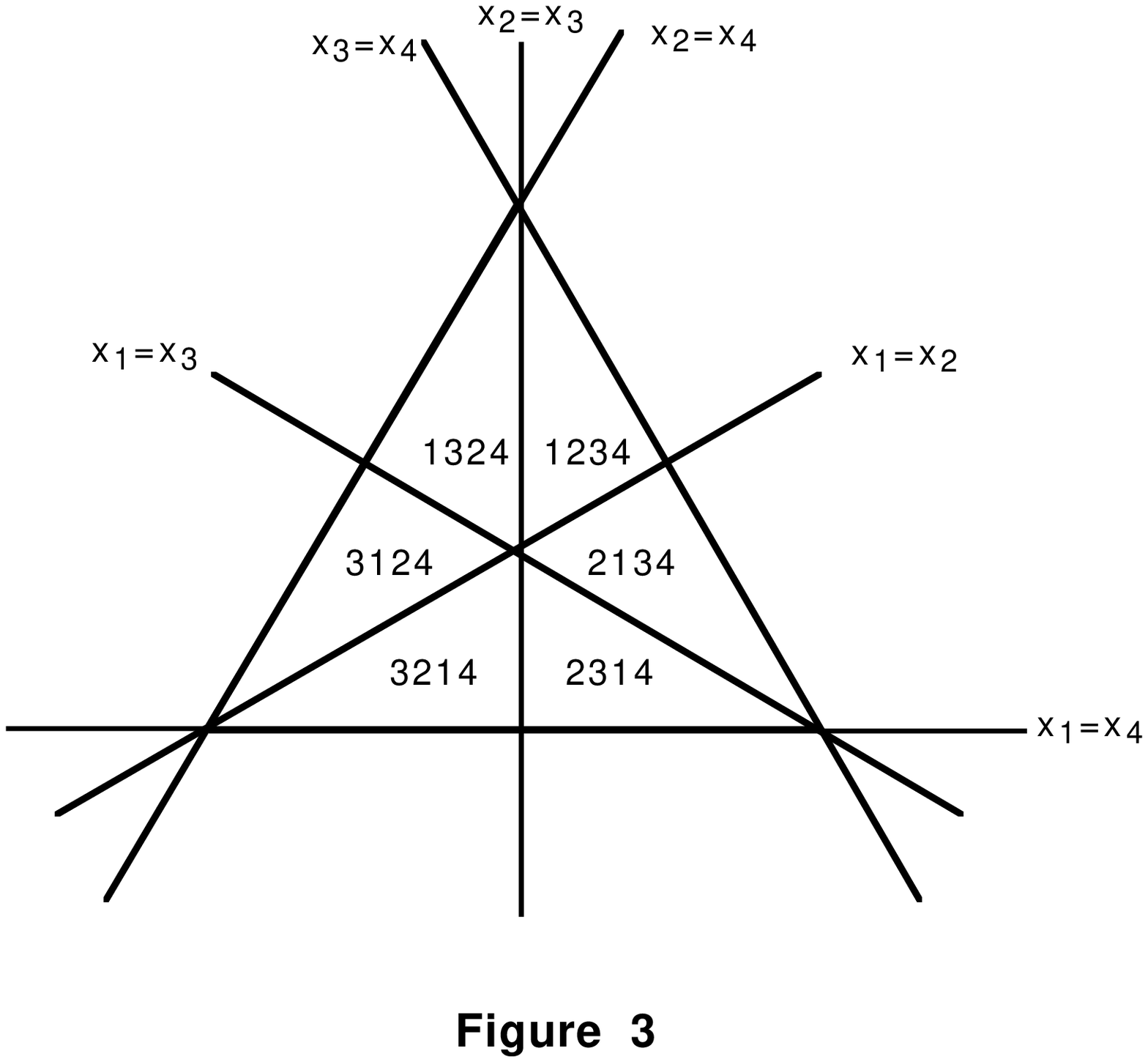}
\end{center}

For type B, when  the generic affine
hyperplane
$H$ is chosen appropriately, one gets the type $B$ splitting basis consisting of  cycles 
$\rho_\omega$ indexed by signed permutations $\omega$ whose right-to-left maxima are
unbarred.   The hyperplane arrangement $\mathcal B_3$ intersected with a cube  is shown in Figure
4.  The regions   that have bounded intersection with $H$ are the ones that are labeled.  The
labels are  the signed permutations whose right-to-left maxima are unbarred.

\vspace{.2in}
\begin{center}
\includegraphics[width=9cm]{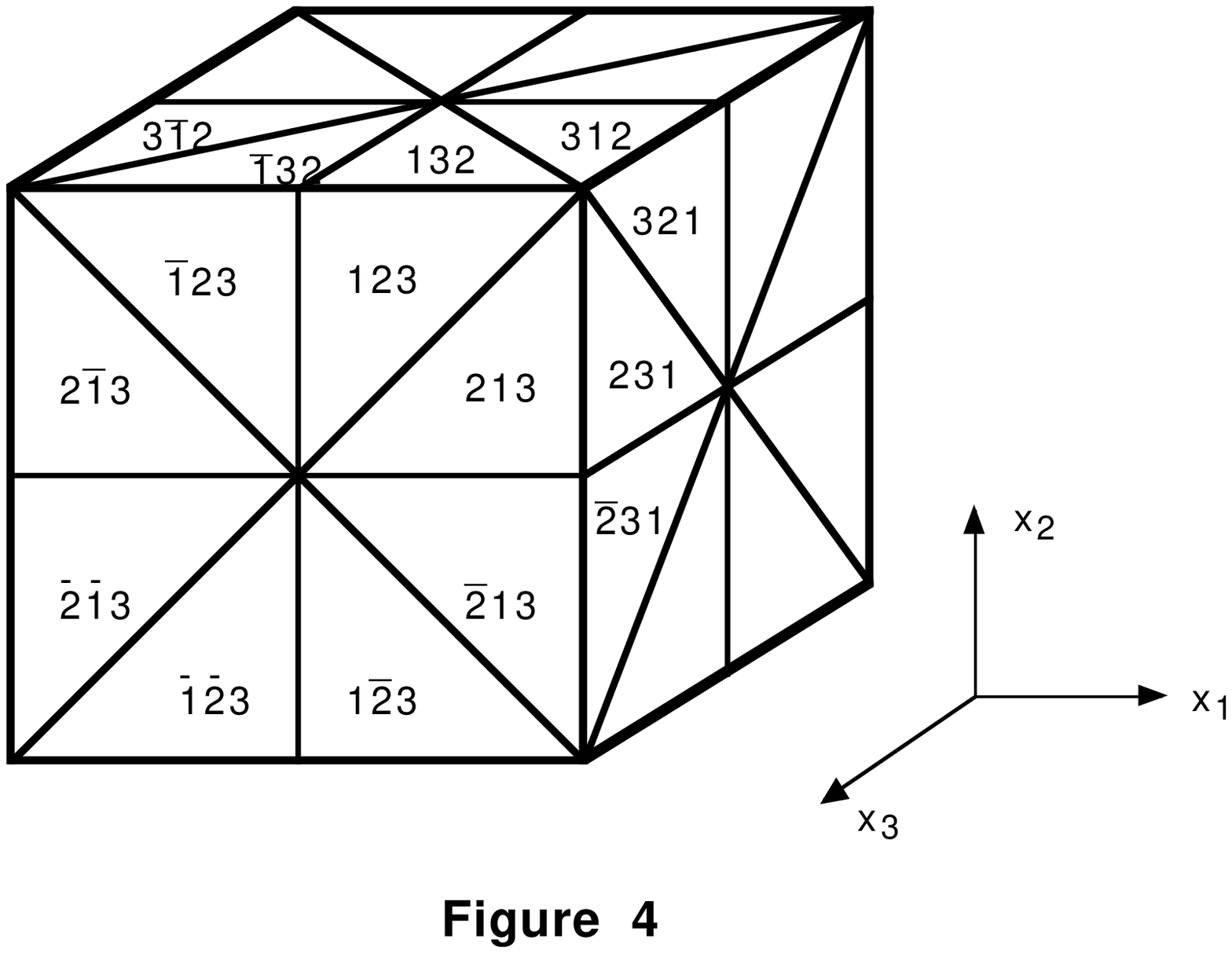}
\end{center}

All arguments in the paper are combinatorial in nature,
which means that they can be carried out for oriented matroids. So the
construction of bases is applicable to geometric lattices of
orientable matroids. Geometrically this means that we can allow some
topological deformation of the hyperplane arrangements.

Major parts of this work (Sections~\ref{affine}, \ref{cent} and~\ref{braid}) were
carried out at the Hebrew University in 1993 during the Jerusalem Combinatorics
Conference.  The rest was added in 1998.
It has been brought to our attention that some of the material in Sections~\ref{affine}
and~\ref{cent} shows similarities with work of others (see e.g. Proposition 5.6 of Damon
\cite{Da96} and parts of Ziegler \cite{Zi88}, \cite{Zi92}); however, there is no
substantial overlap or direct duplication.

\section{A lemma on shellable posets} 

The concept of a {\em shellable complex\,} and a {\em shellable poset\,}
will be considered known. See \cite{BW96} for the definition and basic
properties. In particular, we will make use of the {\em shelling
basis\,} for homology and cohomology \cite[Section~4]{BW96}. A facet
$F$ will be called a {\em full restriction facet\,} with respect to a
shelling if $\RR(F)=F$, where $\RR(\cdot)$ is the restriction operator
induced by the shelling. ({\em Remark:} Such facets were called {\em
homology facets\,} in \cite[Section~4]{BW96}.)

Our notation for posets is that of \cite[Section~5]{BW96}. For
instance, if $P$ is a bounded poset with top element $\hat 1$ and
bottom element $\hat 0$ then $\oli P$ denotes the {\em proper part} of $P$,
 which is defined to be $P\ssm\{\hat 0,
\hat 1\}$; and if
$P$ is an arbitrary poset then $\wh{P}=P\uplus \{\hat 0, \hat 1\}$.
Also, define $P_{<x}:=\{y\in P | y<x\}$ and $P_{\le x}:=\{y\in P | y\le x\}$.

The following simple lemma is a useful devise for identifying bases for homology of simplicial
complexes.  It is used implicitly in
\cite[proof of Theorem~2.2]{Wa96} and variations of it are used in \cite{Br, BrW, GW}. For any
element
$\rho$ of the chain complex of a simplicial complex $\Delta$ and face $F$ of $\Delta$, we denote the
coefficient of  $F$ in $\rho$  by $\langle\rho,F\rangle$. 

\begin{Lem}\label{L1}  Let $\Delta$ be a $d$-dimensional simplicial complex for which $\wt
H_d(\Delta) $ has rank $t$.  If  $\rho_1,\rho_2,\dots,\rho_t$ are $d$-cycles and  $F_1, F_2,
\dots, F_t$ are facets such that the matrix $(\langle \rho_i,F_j\rangle)_{i,j \in [t]}$ is invertible
over
$\Z$, then $\rho_1,\rho_2,\dots,\rho_t$ is a basis for $\wt
H_d(\Delta) $.
\end{Lem}

\begin{proof} Let $\sum_{i=1}^t a_i \rho_i
= 0$.  Then $$(a_1,\dots,a_t) (\langle\rho_i,F_j\rangle)_{i,j \in [t]} = (0,\dots,0).$$
Since $(\langle\rho_i,F_j\rangle)_{i,j \in [t]}$ is invertible, $a_i = 0$ for all $i$.  Hence
$\rho_1,\rho_2,\dots,\rho_t$ are independent over $\Q$ as well as $\Z$.  It follows that
$\rho_1,\rho_2,\dots,\rho_t$ forms a basis over $\Q$.

To see that $\rho_1,\rho_2,\dots,\rho_t$ spans $\wt
H_d(\Delta) $, let $\rho$ be a $d$-cycle.  Then $\rho = \sum_{i=1}^t c_i \rho_i$ where $c_i \in
\Q$.  We have $$(c_1,\dots,c_t) (\langle\rho_i,F_j\rangle )_{i,j \in [t]} =
(\langle\rho,F_1\rangle,\dots,\langle\rho,F_t\rangle)$$ It follows that $$(c_1,\dots,c_t)  =
(\langle\rho,F_1\rangle,\dots,\langle\rho,F_t\rangle)(\langle\rho_i,F_j\rangle)_{i,j \in [t]}^{-1}
\in \Z^t.$$  Hence $\rho$ is in the $\Z$-span of $\rho_1,\rho_2,\dots,\rho_t$.
\end{proof}

Suppose that $\Om$ is a shelling order of the maximal chains of a pure shellable poset
$P$  of length $r$. Let $M$ be the set of maximal elements of $P$.  Recall the following two facts:

\begin{enumerate}
\item[(i)]{For each $m \in M$, a shelling order $\Om^{<m}$ is induced on the maximal chains of
$P_{<m}$ by restricting $\Om$ to the chains containing $m$
\cite[Prop 4.2]{Bj80}.}
\item[(ii)]{A shelling order $\Om^{P\ssm M}$ is induced on the maximal
chains of $P\ssm M$ as follows. Map each maximal chain $c$ in $P\ssm
M$ to its $\Om$-earliest extension $\vap(c)=c \cup \{m\}$, $m\in
M$. Note that $\vap$ is injective. Now say that $c$ precedes $c'$ in
$\Om^{P\ssm M}$ if and only if $\vap(c)$ precedes
$\vap(c')$ \cite[Th. 4.1]{Bj80}.}
\end{enumerate}
Let $\FF(P_{<m})$ and $\FF(P\ssm M)$ denote the sets of full
restriction facets induced by $\Om^{<m}$ and $\Om^{P\ssm M}$. 
Recall from
\cite[Section 4]{BW96} that the shelling $\Omega^{<m}$ induces a basis $B(P_{<m}):=\{\rho_F\}_{F\in
\FF(P_{<m})}$ of $\wt H_{r}(P_{<m})$ which is characterized by the property that $\langle \rho_F,
F^\prime\rangle =
\delta_{F,F^\prime}$ for all $F,F^\prime \in \FF(P_{<m})$.

\begin{Lem}\label{L2}

Let $P$ be a pure poset of length $r$ and $M$ the set of its maximal
elements. Suppose that $P$ is shellable and acyclic.  Then

\begin{enumerate}
\item[(i)] $\FF(P\ssm M) = \biguplus_{m\in M} \FF (P_{<m})$,
\item[(ii)] $\biguplus_{m\in M} B(P_{<m})$ is a basis for 
$\wt{H}_{r-1} (P\ssm M)$.
\end{enumerate}

\end{Lem}

\begin{proof}[Proof of (i)]

 We claim that 
 \beq\label{claim}
c \in \FF(P_{<m}) \implies \varphi(c) = c\cup\{m\} \text{ and } c \in \FF(P \ssm M).\eeq
Let $c\in\FF(P_{<m})$. This means that $c\ssm\{x\}$ is contained in an
$\Om^{<m}$-earlier maximal chain of $P_{<m}$, for every $x\in c$. If
$\vap (c)=c\cup\{m'\}$ with $m'\neq m$ then it would follow that
$c\cup\{m\}$ is a full restriction facet of $P$, contradicting the
assumption that $P$ is acyclic. Hence $\vap(c)=c\cup\{m\}$.  We can also conclude that
$c \in \FF(P\ssm M)$.

It follows from (\ref{claim}) that the sets  $\FF(P_{<m}),\, m \in M$, are disjoint and that
$$
\FF(P\ssm M)\spst \biguplus_{m\in M} \FF(P_{<m}).
$$
The reverse inclusion will be a consequence of the following
computations using the M\"obius function $\mu (\hat0,x)$ of
$\wh{P}$. Since $P$ is acyclic we have that
$$
\sum_{x\in\wh{P}\ssm\{\hat{1}\}} \mu (\hat 0, x)
=
- \mu (\hat 0, \hat 1) = - \wt\chi (P)=0.
$$ 
Hence,
\begin{align*}
|\FF(P\ssm M)|
&=
(-1)^{r} \sum_{x \in \wh P \ssm \{\hat 1 \} \ssm M} \mu (\hat 0, x)\\
&=
(-1)^{r-1} \sum_{m\in M} \mu (\hat 0, m)
=
\sum_{m\in M} |\FF(P_{<m})|.
\end{align*}

\end{proof}

\begin{proof}[Proof of (ii)]

For the homology basis of $\wt{H}_{r}(P\ssm M)$ we will use
Lemma~\ref{L1}. Order $\FF(P\ssm M)$ by $\Om^{P\ssm M}$, and for each
$c\in\FF(P\ssm M) = \uplus_{m\in M} \FF(P_{<m})$, let $m_c$ be defined by
$\vap(c)=c\cup\{m_c\}$. By (\ref{claim}), $c \in \FF(P_{<m_c})$. Let $\rho_c$ be the 
element of $B(P_{<m_c})$ corresponding to $c$.
So, $\rho_c$ is the $(r-1)$-cycle in $P_{<m_c}$ with
coefficient $+1$ at $c$ and coefficient 0 at all
$c'\in\FF(P_{<m_c})\ssm\{c\}$. 

Suppose that $\rho_c$ has nonzero coefficient at some chain
$c'\neq c$. Since $c'$ must come before $c$ in $\Om^{<m}$ (the cycle
$\rho_c$ has support on a subset of the chains in $P_{<m}$ that were
present at the stage during the shelling $\Om^{<m}$ when $c$ was
introduced), it follows that $\vap(c')$ precedes $\vap(c)$ in $\Om$,
and hence that $c'$ precedes $c$ in $\Om^{P\ssm M}$.  
Hence the matrix $(\langle \rho_c, c^\prime \rangle)_{c,c^\prime \in \FF(P\ssm M)}$ is
lower triangular with $1$'s on the diagonal.
It now follows from  Lemma \ref{L1} that
$\biguplus_{m \in M} B(P_{<m})  = \{\rho_c\}_{c\in\FF(P\ssm M)} $ is a basis for $\wt{H}_{r}$
$(P\ssm M)$.
\end{proof}

\section{Affine hyperplane arrangements} \label{affine}

Let $\Aa=\{H_1, \dots, H_t\}$ be an arrangement of affine (or linear) hyperplanes
in $\R^d$. 
Each hyperplane $H_i$ divides $\R^d$ into three components: $H_i$
itself and the two connected components of $\R^d\ssm H_i$. For $x,
y\in\R^d$, say that $x\equiv y$ if $x$ and $y$ are in the same
component with respect to $H_i$, for all $i=1, \dots, t$. This
equivalence relation partitions $\R^d$ into open cells.

Let $P_\Aa$ denote the poset of cells (equivalence classes under
$\equiv$), ordered by inclusion of their closures. $P_\Aa$ is called
the {\em face poset\,} of $\Aa$. It is a finite pure poset with at most $d+1$
rank levels corresponding to the dimensions of the cells. The maximal
elements of $P_\Aa$ are the regions of $\R^d\ssm\bigcup\Aa$. 
See Ziegler \cite{Zi88} for a detailed discussion of these
facts.

Assume in what follows  that the face poset $P_\Aa$ has length $d$.  We will make use of the
following technical properties of 
the order complex of
$P_\Aa$.

\begin{Prop}[{\cite[Section~3]{Zi88}}]\label{P1}\quad
\begin{enumerate}
\item[(i)]{$P_\Aa$ is shellable.}
\item[(ii)]{$P_\Aa$ is homeomorphic to the $d$-ball.}
\item[(iii)]{Let $R$ be a region of $\R^d\ssm\bigcup\Aa$. Then
\begin{equation*}
(P_\Aa)_{<R} \cong
\begin{cases}
(d-1)\text{-sphere}
&\,\text{if $R$ is bounded}\\
(d-1)\text{-ball}
&\,\text{otherwise.}
\end{cases}
\end{equation*}}
\end{enumerate}
\end{Prop}

If $R$ is a bounded region then its closure $\cl(R)$ is a convex
$d$-polytope, and the open interval $(P_\Aa)_{<R}$ is the proper part
of the face lattice of $\cl(R)$. The order complex of $(P_\Aa)_{<R}$,
being a simplicial $(d-1)$-sphere, supports a unique (up to sign) {\em fundamental
$(d-1)$-cycle} $\tau_R$.

Let $\bbarp_\Aa=\{\sigma\in P_\Aa\mid\dim\sigma<d\}$. Equivalently,
$\bbarp_\Aa$ is the poset $P_\Aa$ with its maximal elements
(the regions) removed. Also, let
$\BB=\{\text{bounded regions}\}$.

\begin{Prop}\label{P2}\quad
\begin{enumerate}
\item[(i)]{$\bbarp_\Aa$ has the homotopy type of a wedge of
$(d-1)$-spheres.}
\item[(ii)]{$\{\tau_R\}_{R\in\BB}$ is a basis for $\wt{H}_{d-1}(\bbarp_\Aa)$.}
\end{enumerate}
\end{Prop}

\begin{proof}

Part (i) follows from the fact that shellability is preserved by
rank-selection \cite[Th. 4.1]{Bj80}, and that a shellable pure
$(d-1)$-complex has the stated homotopy type. Since $\{\tau_R\}$ is
(due to uniqueness) the shelling basis for $\wt{H}_{d-1} 
((P_\Aa)_{<R})$ when $R\in\BB$, and $\wt{H}_{d-1} ((P_\Aa)_{<R})=0$
when $R\not\in\BB$, part (ii) follows from Lemma~\ref{L2}.
\end{proof}

Let $L_\Aa$ denote the {\em intersection semilattice\,} of $\Aa$. Its
elements are the nonempty intersections $\bigcap\Aa'$ of subfamilies
$\Aa'\sbst\Aa$, and the order relation is reverse inclusion. $L_\Aa$
is a pure poset of length $d$. Its unique minimal element
is $\R^d$ (corresponding to $\Aa'=\varnothing$), which (according to
convention) will be denoted by $\hat 0$. The minimal elements of
$L_\Aa\ssm\{\hat0\}$ are the hyperplanes $H_i\in\Aa$, and the maximal
elements are the single points of $\R^d$ obtainable as intersections
of subfamilies $\Aa'\sbst\Aa$. $L_\Aa$ is a {\em geometric
semilattice\,} in the sense of \cite{WW86}.

For each cell $\sigma\in P_\Aa$, let $z(\sigma)$ be the affine span of
$\sigma$. The subspace $z(\sigma)$ can also be described as
follows. By definition, $\sigma$ is the intersection of certain
hyperplanes in $\Aa$ (call the set of these hyperplanes $\Aa_\sigma$)
and certain halfspaces determined by other hyperplanes in $\Aa$. Then,
$z(\sigma)=\bigcap\Aa_\sigma$. This shows that $\dim\sigma=\dim
z(\sigma)$ and that $z(\sigma)\in L_\Aa$. The map
$$
z:P_\Aa\to L_\Aa
$$
is clearly order-reversing, and it restricts to an order-reversing map
$$
z:\bbarp_\Aa\to L_\Aa\ssm\{\hat 0\}.
$$

In various versions, the following result appears in several places in the
literature; see the discussion following Lemma~3.2 of \cite{Zi92}.

\begin{Prop}\label{P3}

The map $z:\bbarp_\Aa \to L_\Aa\ssm\{\hat 0\}$ induces homotopy equivalence of
order complexes.
\end{Prop}

\begin{proof}
We will use the Quillen fiber lemma \cite{Q78}. This reduces the
question to checking that every fiber $z^{-1}((L_\Aa)_{\ge x})$ is
contractible, $x\in L_\Aa\ssm\{\hat 0\}$. But by Proposition \ref{P1} (ii) such
a fiber is homeomorphic to a $\dim (x)$-ball, so we are done.
\end{proof}

The simplicial map $z$ induces a homomorphism
$$
z_*:\wt{H}_{d-1} (\bbarp_\Aa)\to \wt{H}_{d-1} (L_\Aa\ssm\{\hat 0\}),
$$
which (as a consequence of Proposition \ref{P3}) is an isomorphism. 
The
following is an immediate consequence of Propositions \ref{P2} and \ref{P3}.

\begin{Thm}\label{T1}

$\{z_* (\tau_R)\}_{R\in\BB}$ is a basis of $\wt{H}_{d-1}
(L_\Aa\ssm\{\hat0\})$.
\end{Thm}

Recall that  $\tau_R$ is the fundamental cycle of the proper part of the face
lattice of the convex polytope $\cl(R)$, for each bounded region
$R$. Since the map $z$ is injective on each lower interval
$(P_\Aa)_{< R}$ it follows that the cycles $z_* (\tau_R)$ are also
``polytopal'', arising from copies of  the proper part of the dual face lattice of $\cl(R)$
embedded in $L_\Aa$. 

\begin{Rem}
{\rm It is a consequence of Theorem \ref{T1} that $$\rank \;\wt{H}_{d-1} 
(L_\Aa\ssm\{\hat0\})=\text{card }\BB.$$ This enumerative corollary is equivalent
to the following result of Zaslavsky \cite{Za75}:
$$
\text{card }\BB=\Big|\sum_{x\in L_\Aa} \mu\(\hat0,x\)\Big|.
$$
Indeed, we have that
$$
\rank \; \wt{H}_{d-1} \(L_\Aa\ssm\{\hat0\}\)
=
\big|\mu_{L_\Aa\cup\{\hat1\}} \(\hat0,\hat1\)\big|,
$$
since $L_\Aa\cup\{\hat1\}$ is the intersection lattice of a central arrangement and is hence a
geometric lattice.  Since
$$
\mu_{L_\Aa\cup\{\hat1\}} \(\hat0,\hat1\)
=
-\sum_{x\in L_\Aa} \mu \(\hat0,x\),
$$
the results are equivalent.}
\end{Rem}

\begin{Rem}
{\rm Our work in this section has the purpose to provide
a short but exact route to the results of the following section,
in particular to Theorem~\ref{T2}. In the process,  a natural
method for constructing bases for geometric semilattices that are
intersection lattices of real affine hyperplane arrangements is given
by Theorem~\ref{T1}.  
For general geometric semilattices, a method for constructing bases  which generalizes
the broken circuit construction of \cite{Bj82} is  given by Ziegler
\cite{Zi92}.  This construction does not reduce to the construction given by 
Theorem~\ref{T1} in the case that the geometric semilattice is the intersection
lattice of a real affine hyperplane arrangement.}
\end{Rem}

\section{Central hyperplane arrangements} \label{cent}

Let $\Aa$ be an essential arrangement of linear hyperplanes in $\R^d$. 
As before, let $L_\Aa$ denote the set of
intersections $\bigcap\Aa'$ of subfamilies $\Aa'\sbst\Aa$ (such
intersections are necessarily nonempty in this case) partially
ordered by reverse inclusion. The finite lattice $L_\Aa$ is called the
{\em intersection lattice\,} of $\Aa$. It is a geometric lattice of
length $d$.

Now, let $H$ be an affine hyperplane in $\R^d$ which is generic with
respect to $\Aa$. Genericity here means that $\dim(H\cap X)=\dim(X)-1$
for all $X\in L_\Aa$. Equivalently, $0\not\in H$ and $H\cap
X\neq \varnothing$ for all $1$-dimensional subspaces $X\in L_\Aa$.

Let $\Aa_H=\{H\cap K\mid K\in\Aa\}$. This is an affine hyperplane
arrangement induced in $H\cong\R^{d-1}$. We denote by $L_{\Aa_H}$ its
intersection semilattice.

\begin{Lem}\label{L3}

$L_{\Aa_H}\cong L_\Aa \ssm\{\hat 1\}$.
\end{Lem}

\begin{proof}

The top element $\hat 1$ of $L_\Aa$ is the 0-dimensional
subspace $\{0\}$ of $\R^d$. Thus $X\mapsto H\cap X$ defines an
order-preserving map $L_\Aa \ssm\{\hat 1\}\to L_{\Aa_H}$, which is
easily seen to be an isomorphism.
\end{proof}

The connected components of $\R^d\ssm\cup\Aa$ are pointed open convex
polyhedral cones, that we call {\em regions}. Although none of these regions are bounded (since $\Aa$
is central), each region
$R$, nevertheless, induces a cycle $\rho_R$ in $\widetilde{H}_{d-2} (\oli L_\Aa)$ as follows. 
Let
$P_R$ denote the face lattice of the closed cone $\cl(R)$. That is,  $P_R$ is the lower interval
 $(P_\Aa)_{\le R}$.
 Clearly $P_R$ is isomorphic to the face lattice of the  convex polytope 
$\cl(R\cap M)$, where $M$ is any
affine hyperplane such that $R\cap M$ is nonempty and bounded.
The map $z:P_\Aa \to L_\Aa$ defined in
Section~\ref{affine} clearly embeds a copy of the dual of  $P_R$ in $L_\Aa$.  Hence
the image
$z(P_R)$ is a subposet of
$L_\Aa$ whose proper part is $(d-2)$-spherical (meaning that its order complex is
homeomorphic to $S^{d-2}$). Let $\rho_R$ be the fundamental cycle (uniquely defined up to sign) of
the proper part of the  subposet $z(P_R)$.

\begin{Thm}\label{T2}
Let $\Aa$ be a central and essential hyperplane arrangement in $\R^d$ and let $H$ be an affine
hyperplane, generic with respect to $\Aa$. Then the collection of cycles $\rho_R$ corresponding
to regions $R$ such that $R\cap H$ is nonempty and bounded, form a basis of
$\wt{H}_{d-2} (\oli{L}_\Aa)$.
\end{Thm}

\begin{proof}

This follows immediately from Theorem \ref{T1} and the fact that
$\oli{L}_\Aa\cong L_{\Aa_H} \ssm\{\hat 0\}$ (Lemma \ref{L3}).
\end{proof}

In order to apply Theorem~\ref{T2} to the examples given in subsequent sections we will need to choose
an appropriate  generic affine hyperplane and determine the regions whose affine slices
are bounded.  The following lemma provides a useful way of doing this.

\begin{Lem} \label{dot} Let $\Aa$ be a central and essential hyperplane arrangement in $\R^d$.
Suppose $\bold v$ is a nonzero element of $\R^d$ such that the affine hyperplane $H_{\bold v}$
through
$\bold v$ and  normal to
$\bold v$, is generic with respect to $\Aa$.  Then for any region $R$ of $\Aa$, $R \cap  H_{\bold v}$
is nonempty and bounded if and only if $\bold v \cdot \bold x > 0$ for all $\bold x  \in R$.
\end{Lem} 

\begin{proof} $(\Rightarrow)$
Suppose
$R
\cap  H_{\bold v}$ is nonempty  and bounded.
It is not difficult to see that if an affine slice of a cone is nonempty and bounded, then the cone
is a cone over the affine slice.  Hence
$R$ is a cone over 
$R
\cap H_{\bold v}$.  That is, every element of $R$ is a positive scalar multiple of an element of
$R \cap H_{\bold v}$.  It follows that since $\bold v \cdot \bold x > 0$ for all $\bold x \in
H_{\bold v}$,
$\bold v
\cdot \bold x > 0$ for all $\bold x \in R$.

$(\Leftarrow)$  Suppose $R\cap H_{\bold v}$ is empty or unbounded.  If the former holds then $\bold v
\cdot \bold x
\le 0$ for all $\bold x \in R$.  Indeed, if $\bold v \cdot \bold x > 0$ for some $\bold x \in R$
then
$\frac{\bold v\cdot \bold v}{\bold v
\cdot \bold x} \bold x
\in R\cap H_{\bold v}$.  

We now assume $R \cap H_{\bold v}$ is unbounded.  Then there is a sequence of points
$\bold x_1,\bold x_2\dots$ in 
$R
\cap H_{\bold v}$ whose distance from the origin goes to infinity.  Let $\bold e_i$ be the  unit
vector in the direction of the vector $\bold x_i$.  Each $\bold e_i$ is in the intersection of R and
the unit sphere centered at the origin.  Hence, by passing to a subsequence if necessary, we can
assume that the sequence of $\bold e_i$'s converge to a unit vector $\bold e$ in the closure of $R$.
Since the cosine of the angle between $\bold e_i$ and
$\bold v$ is
$\frac{\|\bold v\|} { \|\bold x_i\|}$, the  cosine of the angles approach $0$.  Hence the cosine of
the angle between $\bold e$ and $\bold v$ is $0$, or equivalently $\bold v \cdot \bold e = 0$. 

Since $\bold e$ is in the closure of $R$, either $\bold e \in R$ or there is a unique face $F$  of $R$
such that
$\bold e$ is in the interior of $F$. If $\bold e \in R$ we are done.  So suppose $\bold e$ is in the
interior of the face $F$.  If $\bold v \cdot \bold x = 0$ for all $\bold x \in F$ then the linear span
of
$F$ is an intersection of hyperplanes contained in the linear hyperplane with normal vector $\bold
v$.  This contradicts the genericity of $H_{\bold v}$.  It follows that $\bold v \cdot \bold x \ne 0$
for some
$\bold x
\in F$.  If
$\bold v
\cdot \bold x < 0$ then there is a point $\bold  y \in R$ that is close enough to $\bold x$ so that
$\bold v
\cdot \bold y < 0$ and we are done.  If
$\bold v
\cdot \bold x > 0$ then consider the point $\bold e - a \bold x$ where $a > 0$.  We have $\bold v
\cdot (\bold e-a \bold x) = -a (\bold  v
\cdot \bold x) < 0$.  By choosing $a$ to be small enough, we insure that the point $\bold e-a \bold x$
is close enough to $\bold e$ to be in 
$F$, since $\bold e$ is in the interior of $F$.  Hence we have a point in $F$ whose dot product with
$\bold v$ is negative, putting us back in the previous case.
\end{proof}

\begin{Rem} {\rm Theorem~\ref{T2} can be extended to a geometric
construction of bases for the Whitney homology, or equivalently
the Orlik-Solomon algebra,  of the intersection lattice of a real central hyperplane
arrangement. This involves the definition of a vector $v$
being {\em totally generic} with respect to the arrangement.
Since we will not pursue this direction we omit further
mention of it.}
\end{Rem}

\section{Oriented matroids} 

The arguments and results of the previous two sections can be
generalized to oriented matroids. This generalization will be outlined
in this section. The treatment here will be sketchy and can be
skipped with no loss of continuity. The basics of oriented matroid
theory will be assumed to be known. We refer to \cite{OM} for all
definitions and notation.

Let $(\LL, E, g)$ be an {\em affine oriented matroid\,} of rank $r$
and with {\em affine face lattice\,} $\LL^+=\{X\in\LL\mid X_g=+\}$,
cf. \cite[Section~4.5]{OM}. The maximal elements of $\LL^+$ are the
{\em topes}, corresponding to regions in the realizable case. Let
$\LL^{++}$ be the {\em bounded complex} (a subcomplex of $\LL^+$), and
let $\BB^{++}$ be the set of bounded topes, i.e.,
$\BB^{++}=\{X\in\LL^{++}\mid\rank(X)=r-1\}$.

We have from \cite[Th.~4.5.7]{OM} that $\LL^+$ is a shellable
ball. Furthermore, if $T\in\BB^{++}$ then the order complex of the
open interval $(0, T)$ in $\LL^+$ is homeomorphic to
$S^{r-2}$ \cite[Cor.~4.3.7]{OM}. 
Therefore, each $T\in\BB^{++}$ induces a spherical
fundamental cycle $\tau_T$ in $\wt{H}_{r-2}({\bbarl})$, where $\bbarl
= \LL^+\ssm\{\text{topes}\}$. 

\begin{Prop}\label{P4}\quad
\begin{enumerate}
\item[(i)]{$\bbarl$ has the homotopy type of a wedge of $|\BB^{++}|$
copies of the $(r-2)$-sphere.}
\item[(ii)]{$\{\tau_T\}_{T\in\BB^{++}}$ is a basis of $\wt{H}_{r-2}
(\bbarl)$.}
\end{enumerate}
\end{Prop}

\begin{proof}
The proof of Proposition \ref{P2} generalizes.
\end{proof}

Now, let $L$ be the {\em intersection lattice\,} (or ``lattice of
flats'') of the oriented matroid $\LL$, and let $z:\LL\to L$ be the
``zero map'' \cite[Prop. 4.1.13]{OM}. Furthermore, let $L^{g}:=\{x\in
L\mid g \not\in x\}=L\ssm[g,\hat 1]$. This is a geometric
semilattice. The zero map restricts to an order-reversing surjection
$z:\LL^+\to L^{g}$, and further to a surjection $z:\bbarl \to L^{g}\ssm\{\hat 0\}$.

\begin{Prop}\label{P5}

The map $z:\bbarl\to L^{g}\ssm\{\hat 0\}$ induces homotopy
equivalence of order complexes.
\end{Prop}

\begin{proof}
The proof of Proposition \ref{P3} generalizes. Here one uses that the
Quillen fibers $z^{-1} ((L^{g})_{\ge x})$, $x\neq \hat 0$, are balls
by \cite[Th.~4.5.7]{OM}, and hence contractible.
\end{proof}

The restriction of $z$ to an open interval $(0,T)$ in $\LL^+$, with
$T\in\BB^{++}$, gives an isomorphism of $(0,T)$ onto its image in
$L^{g}$. This image is a subposet of $L^{g}\ssm\{\hat 0\}$
homeomorphic to the $(r-2)$-sphere. Let $\rho_T\in\wt{H}_{r-2}(L^{g}\ssm\{\hat 0\})$ be the corresponding fundamental cycle.

\begin{Thm}\label{T3} \quad
\begin{enumerate}
\item[(i)]{$L^{g}\ssm\{\hat 0\}$ has the homotopy type of
$|\BB^{++}|$ copies of the $(r-2)$-sphere.}
\item[(ii)]{$\{\rho_T\}_{T\in\BB^{++}}$ is a basis of $\wt{H}_{r-2}(L^{g}\ssm\{\hat 0\})$.}
\end{enumerate}
\end{Thm}

\begin{proof}

This follows from Propositions \ref{P4} and \ref{P5}, since $$z_*:\wt{H}_{r-2}
(\bbarl)\to\wt{H}_{r-2} (L^{g}\ssm\{\hat 0\})$$ is an isomorphism and
$z_*(\tau_T)=\rho_T$.
\end{proof}

The treatment of affine oriented matroids so far parallels that of
affine hyperplane arrangements in Section~3. We will now move on to
the oriented matroid version of the material in Section~4.

Let $\LL\sbst\{+, -, 0\}^E$ be an oriented matroid of rank $r$, and
let $z:\LL\to L$ be the zero map to the corresponding intersection
lattice $L$. Let $\LL_g\sbst\{+,-,0\}^{E\uplus g}$ be an extension of
$\LL$ by a generic element $g\not\in E$. Genericity here means that
$g\not\in\text{\,span\,}A$ for every $A\subseteq E$ with $\rank(A)<r$,
cf. \cite[Sect. 7.1]{OM}.

Consider the affine oriented matroid $(\LL_g, E\uplus g, g)$ and let
$L_g$ be its intersection semilattice. We have that $z(\LL_g)=L_g$ and
$z(\LL^+_g)=(L_g)^g$.

\begin{Lem}\label{L4}

$(L_g)^g\cong L\ssm\{\hat 1\}$.
\end{Lem}

\begin{proof}

This analog of Lemma \ref{L3} is clear. It is basically a reformulation of
the definition of genericity.
\end{proof}

Let $\BB^{++}$ be the bounded topes of $\LL_g$ (with respect to
$g$). Because of the isomorphism
$$
(L_g)^g\ssm\{\hat 0\}
\cong
\oli{L}:=L\ssm\{\hat 0, \hat 1\},
$$
we get cycles $\rho_T \in \wt{H}_{r-2}(\oli{L})$ as before.

\begin{Thm}\label{T4}\quad
\begin{enumerate}
\item[(i)]{$\oli{L}$ has the homotopy type of a wedge of
$|\BB^{++}|$ copies of the $(r-2)$-sphere.}
\item[(ii)]{$\{\rho_T\}_{T\in\BB^{++}}$ is a basis for
$\wt{H}_{r-2}(\oli{L})$.}
\end{enumerate}
\end{Thm}

\begin{proof}

This follows from Theorem \ref{T3} and Lemma \ref{L4}.
\end{proof}

The theorem gives a geometric method for constructing a basis for the
homology of the geometric lattice of {\em any\,} orientable
matroid. Note that to define the set $\BB^{++}$, and hence the basis,
we must make a generic extension of $\LL$. Different extensions will
yield different bases.

\section{Type $A$: The braid arrangement} \label{braid}

The  hyperplane arrangement $\Aa_{n-1} = \{H_{ij} : 1 \le i < j \le n\}$  in $\R^n$,
where 
$H_{ij} =
\{\bold x \in
\R^n : x_i =x_j\}$, is known as the {\em braid arrangement} or the {\em type $A$ Coxeter
arrangement}.   The  orthogonal reflection
$\sigma_{ij}$ across  the hyperplane
$H_{ij}$ acts on
$(x_1,\dots,x_n)
\in
\R^n$ by switching its
$i$th and
$j$th coordinates.  These reflections generate the symmetric group $S_n$ acting on $\R^n$ by permuting
coordinates.  

The braid arrangement is not essential.  To make it essential let
$$K = \{\bold x \in \R^n : x_1+\dots+x_n =0 \}$$
and define $$\Aa^\prime_{n-1} = \{H^\prime_{ij} = K \cap H_{ij} : 1 \le i < j \le n\}.$$  Then
$\Aa^\prime_{n-1}$ is an essential central hyperplane arrangement in the $(n-1)$-dimensional space
$K$.  It is clear that the intersection lattices $L_{\Aa_{n-1}}$ and  $L_{\Aa^\prime_{n-1}}$ are
isomorphic. They are also isomorphic to the partition lattice $\Pi_n$.  Indeed, for each $\pi \in
\Pi_n$, let $\ell_\pi$ be the linear subspace of $\R^n$ consisting of all points $(x_1,\dots, x_n)$
such that $x_i = x_j$ whenever $i$ and $j$  are in the same block of $\pi$. The map $\pi \mapsto
\ell_\pi\cap K$ is an isomorphism from $\Pi_n$ to $L_{\Aa^\prime_{n-1}}$. Let $\gamma$ denote the
inverse of this isomorphism.

The arrangement $\Aa^\prime_{n-1}$  has $n!$ regions which are all simplicial cones and are in a
natural one-to-one correspondence with the elements of the associated Coxeter group $S_n$.  Under
this correspondence a permutation
$\omega
\in S_n$ corresponds to the region
$$R_\omega = \{\bold x \in K : x_{\omega(1)} < x_{\omega(2)} < \dots < x_{\omega(n)}\}.$$
Consider the cycle $\rho_{R_\omega} \in \wt H_{n-3}(\oli L_{A^\prime_{n-1}})$ whose
general construction was given in Section~\ref{cent}. We now give a simple explicit description of
the image of
$\rho_{R_\omega}$ in $\wt H_{n-3}(\Pi_n)$ under the isomorphism 
$\gamma$. 

To split a permutation $\omega \in S_n$ at positions $i_1 < \cdots < i_k$ in $[n-1]$ is to form the
partition with $k+1$ blocks, 
\begin{equation} \label{split} \{\omega(1),\dots,\omega(i_1)\},
\{\omega(i_1+1),\dots,\omega(i_2)\},\dots,\{\omega(i_k+1),\dots,\omega(n)\}.\end{equation} 
(To split $\omega$ at the empty set of positions is to form the partition with one block.) Let
$\Pi_\omega$ denote the induced subposet of $\Pi_n$ consisting of all partitions obtained by splitting
the permutation
$\omega$.  Clearly $\Pi_\omega$ is isomorphic to the lattice of subsets of $[n-1]$. Hence 
$\oli\Pi_\omega$ is spherical.   

\begin{Prop}\label{funda}  For all $\omega \in S_n$, the image $\gamma(\rho_{R_\omega})$ is the
fundamental cycle of 
$\oli\Pi_\omega$.
 
\end{Prop}

\begin{proof} Recall the map $z: P_{\Aa_{n-1}^\prime} \to L_{\Aa_{n-1}^\prime}$   that takes
cells to their affine span (defined in Section~\ref{affine}).  We show that
$\gamma$ restricts to an isomorphism from  the subposet
$z(P_{R_\omega})$ of
$L_{A^\prime_{n-1}}$ to $\Pi_\omega$.  Observe that elements of the face lattice  $P_{R_\omega}$ are
sets of the form 
\begin{eqnarray*}\{\bold x\in K : x_{\omega(1)} = \dots = x_{\omega(i_1)} < x_{\omega(i_1+1)} =
\cdots = x_{\omega(i_2)}
 < \cdots \\
\cdots < x_{\omega(i_k+1)} = \cdots = x_{\omega(n)}\},\end{eqnarray*}
where $1 \le i_1 < \dots < i_k \le n-1$.
The  linear span of such a set is the subspace
\begin{eqnarray}\label{space}\{\bold x\in K : x_{\omega(1)} = \cdots = x_{\omega(i_1)},\,\,
x_{\omega(i_1+1)} =
\cdots = x_{\omega(i_2)},\,\,
  \dots  \\ \nonumber
\dots ,\,\, x_{\omega(i_k+1)} = \cdots = x_{\omega(n)}\}.\end{eqnarray}
Hence $z(P_{R_\omega})$ is the poset of subspaces of the form given in (\ref{space})
ordered by reverse inclusion.  Clearly $\gamma$ takes the subspace given in (\ref{space}) to the
partition given in (\ref{split}).
\end{proof}

We now choose a vector $ \bold v $ in $K$ that satisfies the hypothesis of 
Lemma~\ref{dot} and use Lemma~\ref{dot} to describe the permutations $\omega \in S_n$ for which the
regions $R_\omega
\cap H_\bold v$ are bounded.

\begin{Prop}\label{bounda} Let $ \bold v  = (-1,-1,\dots,-1,n-1) \in K$.  Then the affine hyperplane
$H_\bold v \cap K$ is generic  with respect to the arrangement $\Aa^\prime_{n-1}$ of
$K$.  Moreover, for all $\omega \in S_n$, $R_\omega
\cap H_\bold v$ is bounded if and only if $\omega(n)=n$.
\end{Prop}

\begin{proof}
Recall that genericity is equivalent to the condition that for all $1$-dimensional subspaces $X \in
L_{\Aa^\prime_{n-1}}$, $H_\bold v \cap X \ne \emptyset$.  The $1$-dimensional intersections of
hyperplanes in $\Aa^\prime_{n-1}$ have the form $$X=\{\bold x \in K : x_{i_1} = x_{i_2} = \cdots =
x_{i_k},\,\, x_{i_{k+1}} = x_{i_{k+2}} = \cdots = x_{i_{n-1}} = x_n \},$$
where $1 \le k \le n-1$ and $\{\{i_1,i_2,\dots,i_k\},\{i_{k+1},\dots,i_{n-1},n\}\}$ is a partition of
$[n]$.   Let
$\bold x
\in
\R^n$ be defined by
$$x_{i_j} = \begin{cases} -\frac{(n-k)(n-1)} k &j= 1,\dots k \\ \quad n-1 &j=k+1,\dots,n-1
\end{cases}$$ and $ x_n = n-1$.  One can easily check that $\bold  x \in H_\bold v \cap X$.  Hence
$H_\bold v \cap K$ is generic.

To prove that $R_\omega \cap H_\bold v$ is bounded if and only if $\omega(n) = n$, we apply
Lemma~\ref{dot}. Note that for all $\bold x \in K$, $$ \bold v  \cdot \bold x = n x_n.$$
Suppose $\omega(n) = n$. For all $\bold x \in R_\omega$, we have  $x_{\omega(1)} < \cdots < x_{\omega(n)}$
and
$\sum_{i=1}^n x_i = 0$. Hence $x_{\omega(n)} > 0$.  It follows that $ \bold v  \cdot \bold x = n x_{\omega(n)} > 0$ for all $\bold x 
\in R_\omega$.  By Lemma~\ref{dot} $R_\omega \cap H_\bold v$ is bounded.  
Now suppose $\omega(n) \ne n$.  Clearly there exists an $\bold x \in R_\omega$ such that
$x_{\omega(i)} < 0 $ for all $i = 1,\dots,n-1$.  For such an $\bold x$, we have $ \bold v  \cdot \bold
x  = n x_n < 0$. We
conclude by  Lemma~\ref{dot} that $R_\omega \cap H_\bold v$ is not bounded.
\end{proof}

\begin{Thm} [Splitting basis \cite{Wa96}]  \label{splita} For  $\omega \in S_n$, 
let $\rho_\omega$ be
the fundamental cycle of
$\oli \Pi_\omega$. Then
$$\{\rho_\omega : \omega
\in S_n \,\,\, {\rm{and}} \,\,\, \omega(n)=n \} $$
is a basis for  $\wt H_{n-3}(\oli \Pi_n)$.
\end{Thm}

\begin{proof}  This follows from Theorem~\ref{T2} and Propositions~\ref{funda} and~\ref{bounda}.
\end{proof}

\section{Type $B$} \label{secb}

The type $B$ Coxeter arrangement is the  hyperplane arrangement \begin{eqnarray*} \BB_n & = &\{ x_i = 
x_j : 1
\le i < j
\le n\} \, \cup \, \{ x_i =  - x_j : 1 \le i < j
\le n\} \, \cup \\ & &  \{x_i = 0 : 1 \le i \le n \}.\end{eqnarray*}   The orthogonal reflections
across the hyperplanes generate the  hyper\-octahedral group $B_n$.  We will view the elements
of $B_n$ as  signed  permutations, that is, words consisting of $n$ distinct letters from $[n]$ where
any of the letters can have a bar placed above it. It will also be convenient to  express elements of
$B_n$  as pairs
$(\omega,\epsilon)$ where $\omega \in S_n$ and $\epsilon \in \{-1,1\}^n$. If $\epsilon_i =1$ then
$\omega(i)$ does not have a bar over it and if $\epsilon(i) = -1$ then $\omega(i)$ has a bar over it.
For example, the signed permutation $\bar 3 5 4 \bar 2 1$ can be expressed  as
$(35421,-1\,1\,1\,-1\,1).
$ A signed permutation $(\omega,\epsilon)$ maps $(x_1,\dots,x_n) \in \R^n$ to $(\epsilon_1
x_{\omega(1)},\dots, \epsilon_n x_{\omega(n)})$.

The arrangement $\BB_n$ is essential and has $2^n n!$ regions which are all simplicial cones and are in
a natural one-to-one correspondence with the elements of the hyperoctahedral group $B_n$.  Under this
correspondence a signed permutation
$(\omega,\epsilon)
\in B_n$ corresponds to the region
$$R_{\omega,\epsilon} = \{\bold x \in \R^n : 0 < \epsilon_1 x_{\omega(1)} < \epsilon_2
x_{\omega(2)} <
\dots < \epsilon _n x_{\omega(n)}\}.$$

The intersection lattice $L_{\BB_n}$ is isomorphic to the signed partition lattice $ \Pi^B_n$
which is defined as follows.  Let $\pi$ be a partition  of the set $\{0,1,\dots,n\}$.  The block
containing $0$ is called the {\it zero block}.  To {\it bar} an element of a block of  $\pi$ is to
place a bar above the element and to unbar a barred element is to remove the bar.  A signed partition
is a partition of the set
$\{0,1,\dots,n\}$ in which any of the nonminimal elements of  any of the nonzero blocks are barred.
For example,
$$057 \mid 1 \bar 2 9 \mid 3\bar 4 \bar 6 8$$
is a signed partition of $\{0,1,\dots,9\}$. 
It will be convenient to sometimes express a barred letter $\bar a$ of a signed partition as $(a,-1)$
and an unbarred letter as $(a,1)$.  

To  bar a  block
$\B$ in a signed partition is to bar all unbarred elements in
$\B$ and to unbar all barred elements in $\B$. We denote this by $\bar \B$. To unbar a block $\B$
is to unbar all barred elements of $\B$.  We denote this by $\wt \B$. For example,
$$\overline{3
\bar 4
\bar 6 8} =
\bar 3 4 6
\bar 8 \quad {\rm and} \quad \widetilde{3
\bar 4
\bar 6 8} = 3468.$$  Let $ \Pi^B_n$ be the poset of signed partitions of $\{0,1,\dots,n\}$ with
order relation defined by $\pi \le \tau$ if for each block $\B$ of  $\pi$, either $\B$ is
contained in a nonzero block of $\tau$,  $\bar \B$ is
contained in a nonzero block of $\tau$ or $\wt \B$ is contained in the zero block of
$\tau$.  For example
$$057\mid 1 \bar 2 9 \mid 3\bar 4 \bar 6 8 < 057 \mid 1\bar 2 9 3 \bar 4 \bar 6  8,$$ 
$$057\mid 1 \bar 2 9 \mid 3\bar 4 \bar 6 8 < 057 \mid 1\bar 2 9 \bar 3 4 6 \bar 8$$
and 
$$057\mid 1 \bar 2 9 \mid 3\bar 4 \bar 6 8 < 0573468 \mid 1\bar 2 9 .$$
The poset $\Pi^B_n$ is an example of a Dowling lattice \cite{D}.

For each signed partition $\pi \in  \Pi^B_n$, let $\ell_\pi$ be the linear subspace of $\R^n$
consisting of all points $(x_1,x_2,\dots,x_n)$ such that 
\begin{itemize}
\item $x_i = x_j$ whenever $i$ and $j$ are in the
same block of $\pi$ and both are barred or both are unbarred,
\item $x_i = -x_j$ whenever $i$ and $j$ are in the
same block of $\pi$ and one is barred and the other is unbarred,
\item $x_i = 0$ whenever $i$ is in the zero block.
\end{itemize}
The map $\pi \mapsto \ell_\pi$ is an isomorphism from $ \Pi^B_n$ to $L_{{\mathcal B}_n}$.
Let $\gamma$ be the inverse of this isomorphism. 

The cycle $\rho_{R_{\omega,\epsilon}} \in \wt H_{n-2}(\oli L_{{\mathcal B}_n})$ maps under $\gamma$
to the fundamental cycle of a spherical subposet of $ \Pi^B_n$, which we now describe.  To split a
signed permutation
$(\omega,\epsilon)$ at positions
$i_1 < \cdots < i_k$ in $\{0,1,\dots,n-1\}$  is to form the signed partition 
with blocks $\B_0,\B_1,\dots,\B_k$ where $\B_0 = \{0,\omega(1),\dots,\omega(i_1)\}$ (here
$\omega(0) = 0$), and  for all $j = 1,\dots,k$, 
 $$\{(\omega(i_j+1),\epsilon_{i_j+1}), \dots, (\omega(i_{j+1}), \epsilon_{i_{j+1}})\}$$ is either
$\B_j$ or $\bar \B_j$ (here $i_{k+1} = n$).  For example, if we split the signed permutation
$$\bar 3 56 \bar 1 87 \bar 4 2$$ at positions $2,5$ we get the signed partition
$$035 \mid \bar 6 1 \bar 8 \mid 7\bar 4 2.$$

For each signed permutation $(\omega,\epsilon) \in B_n$, let $\Pi_{\omega,\epsilon}$ be the induced
subposet of $\Pi^B_n$ consisting of all signed partitions obtained by splitting the
signed permutation $(\omega,\epsilon)$.  Just as for type A, $\oli\Pi_{\omega,\epsilon}$ is
spherical because  $\Pi_{\omega,\epsilon}$ is isomorphic to the lattice of subsets of $[n]$. 

\begin{Prop} \label{fundb}  For all $(\omega,\epsilon) \in B_n$,
the image $\gamma(\rho_{R_{\omega,\epsilon}})$ is the fundamental cycle of 
$\oli\Pi_{\omega,\epsilon}$.
 
\end{Prop}

\begin{proof}  Let $z: P_{\BB_n} \to L_{\BB_n}$ be the map that takes cells to
their affine span (cf. Section~\ref{affine}).  We show that $\gamma$ restricts to an isomorphism
from  the subposet
$z(P_{R_{\omega,\epsilon}})$ of
$L_{\BB_n}$ to $\Pi_{\omega,\epsilon}$.  Observe that elements of the face lattice
$P_{R_{\omega,\epsilon}}$ are sets of the form 
\begin{eqnarray*}\{\bold x\in \R^n \mid 0= x_{\omega(1)} = \dots = x_{\omega(i_1)} <
\epsilon_{i_1+1} x_{\omega(i_1+1)} =
\dots = \epsilon_{i_2} x_{\omega(i_2)}
 < \dots \\
\dots < \epsilon_{i_k+1} x_{\omega(i_k+1)} = \dots =  \epsilon_n x_{\omega(n)}\},\end{eqnarray*}
where $0 \le i_1 < \dots < i_k \le n-1$.
The  linear span of such a set is the subspace
\begin{eqnarray}\label{space2} \\ \nonumber\{\bold x\in \R^n \mid 0 = x_{\omega(1)} = \dots =
x_{\omega(i_1)},\,\,
\epsilon_{i_1+1} x_{\omega(i_1+1)} =
\dots = \epsilon_{i_2} x_{\omega(i_2)},\,\,
  \dots  \\ \nonumber
\dots ,\,\, \epsilon_{i_k+1} x_{\omega(i_k+1)} = \dots = \epsilon_n x_{\omega(n)}\}.\end{eqnarray}
Hence  $z(P_{R_{\omega,\epsilon}})$ is the poset of subspaces of the form given in (\ref{space2})
ordered by reverse inclusion.  Clearly $\gamma$ takes the subspace given in (\ref{space2}) to the
signed partition obtained by splitting $(\omega,\epsilon)$ at positions $i_1,i_2,\dots,i_k$.
\end{proof}

Let
$(\omega,\epsilon)
\in B_n$. We say that
$\omega(i)$ is a {\it right-to-left maximum} of
$(\omega,\epsilon)$ if $\omega(i) > \omega(j)$ for all $j >i$.

\begin{Prop}\label{bounded} Let $\bold v = (1,2,2^2,\dots,2^{n-1})$.  Then the affine hyperplane
$H_\bold v$ is generic  with respect to the arrangement $\BB_n$.  Moreover, for all
$(\omega,\epsilon)
\in B_n$,
$R_{\omega,\epsilon}
\cap H_\bold v$ is bounded if and only if all right-to-left maxima of $(\omega,\epsilon)$ are
unbarred.
\end{Prop}

\begin{proof}
  The $1$-dimensional intersections of hyperplanes
in $\BB_n$ have the form $$X=\{\bold  x \in \R^n \mid 0 =x_{i_1} = x_{i_2} = \cdots = x_{i_k},\,\,
\epsilon_{k+1}x_{i_{k+1}} = \cdots = \epsilon_n x_{i_n} \}$$
where $0 \le k \le n-1$, $\{\{0,i_1,\dots,i_k\},\{i_{k+1},\dots,i_{n}\}\}$ is a partition of
$\{0,1,\dots,n\}$ and $\epsilon_i = \pm 1$, for $i = k+1,\dots,n$.   Let
$\bold y  \in
\R^n$ be defined by
$$y_{i_j} =  \begin{cases} 0 &j= 1,\dots k \\ \epsilon_j
&j=k+1,\dots,n.
\end{cases}$$
  One can easily check that $\frac{4^n-1} {3 \sum_{t=k+1}^n \epsilon_t2^{i_t-1}}\bold y \in H_\bold v
\cap X$.  Hence $H_\bold v \cap X \ne \emptyset$  for all
$1-$dimensional spaces $X$ in $L_{\BB_n}$, which means that $H_\bold v$ is generic.

Now suppose $(\omega,\epsilon) \in B_n$ and $R_{\omega,\epsilon}$ is bounded.  Assume that $\omega(m)$
is a right-to-left maximum that is barred. We will reach a contradiction of Lemma~\ref{dot} by
producing a vector
$\bold y$
in the closure of
$R_{\omega,\epsilon}$ such that $ \bold v  \cdot \bold y < 0$.  Let
$ \bold y \in \R^n$ be such that
$$ y_{\omega(i)} = \begin{cases} 0 & {\rm if} \,\, 1 \le i < m \\ 
\epsilon_i & {\rm if} \,\, m \le i \le n.
\end{cases} $$
Clearly $ \bold y \in \cl(R_{\omega,\epsilon})$.  We have 
\begin{eqnarray*}\bold v \cdot \bold y &=& \sum_{j=m}^n 2^{\omega(j)-1} \epsilon_j
\\  &\le& -2^{\omega(m)-1} + \sum_{j=m+1}^n 2^{\omega(j)-1}.\end{eqnarray*}
Since $\omega(m) > \omega(j)$ for all $j > m$, we have 
$$\sum_{j=m+1}^n 2^{\omega(j)-1} \le
\sum_{i=1}^{\omega(m)-1} 2^{i-1}\,\, = \,\,2^{\omega(m) -1} -1.$$ 
Combining this with the previous
inequality yields $ \bold v  \cdot \bold y \le -1$.

Now suppose $(\omega,\epsilon) \in B_n$ and all right-to-left maxima of $(\omega,\epsilon)$ are
unbarred.    Let $m_1 < m_2 < \cdots < m_k $ be the positions of the
right-to-left maxima.  
Let $\bold x \in R_{\omega,\epsilon}$. Then
\begin{eqnarray*} \bold  v  \cdot \bold x &=& \sum_{i=1}^n 2^{\omega(i)-1} x_{\omega(i)}\\ & =&
\sum_{j=1}^k
\sum_{i=m_{j-1}+1}^{m_j} 2^{\omega(i)-1} x_{\omega(i)},
\end{eqnarray*}
where $m_0 = 0$.   
We have $$\omega(i) < \omega(m_j)$$
for all $i =
m_{j-1}+1, \dots, m_j-1$, since $\omega(m_j)$ is a right-to-left maximum and $\omega(i)$ is not.
Also, 
$$0<\epsilon_i x_{\omega(i)} < x_{\omega(m_j)}$$ 
for all $i =
m_{j-1}+1, \dots, m_j-1$, since $\omega(m_j)$ is unbarred.  Consequently
\begin{eqnarray*}\sum_{i=m_{j-1}+1}^{m_j} 2^{\omega(i)-1} x_{\omega(i)} &\ge& 
2^{\omega(m_j)-1} x_{\omega(m_j)} - \sum_{\scriptsize\begin{array}{c}
i=m_{j-1}+1
\\
\epsilon_i < 0 \end{array}
}^{m_j-1}  2^{\omega(i) -1} \epsilon_i x_{\omega(i)} \\ &\ge&
2^{\omega(m_j)-1} x_{\omega(m_j)} - \sum_{\scriptsize\begin{array}{c}
i=m_{j-1}+1
\\
\epsilon_i < 0 \end{array}
}^{m_j-1}  2^{\omega(i) -1}  x_{\omega(m_j)} \\ &=&
x_{\omega(m_j)}\,\,\Large(2^{\omega(m_j)-1}  - \sum_{\scriptsize\begin{array}{c}
i=m_{j-1}+1
\\
\epsilon_i < 0 \end{array}
}^{m_j-1}  2^{\omega(i) -1} \Large)\\ &\ge& x_{\omega(m_j)}\,\,\Large(2^{\omega(m_j)-1} -
(2^{\omega(m_j)-1} -1)\Large) \\ &>& 0.
\end{eqnarray*}
Thus we have shown  that $ \bold v  \cdot \bold x > 0$. By Lemma~\ref{dot}, $R_{\omega,\epsilon} \cap
H_\bold v$ is bounded.
\end{proof}

\begin{Thm}[type B splitting basis] \label{splitb}    For  $(\omega,\epsilon) \in B_n$, let
$\rho_{\omega,\epsilon}$ be the fundamental cycle of the spherical poset
$\oli\Pi_{\omega,\epsilon}$. Then
$$\{\rho_{\omega,\epsilon} \mid (\omega,\epsilon)
\in B_n  \text{ and all right-to-left maxima of } (\omega, \epsilon) \text{ are unbarred} 
\}
$$ is a basis for  $\wt H_{n-2}(\overline{\Pi^B_n})$.
\end{Thm}

\begin{proof} This follows from Theorem \ref{T2} and Propositions~\ref{fundb} and~\ref{bounded}. 
\end{proof}

\begin{Cor}[Dowling \cite{D}] \label{Dow} The rank of $\wt H_{n-2}(\oli {\Pi^B_n})$ is $$1 \cdot 3
\cdot 5
\cdots (2n-1).$$
\end{Cor}

\begin{proof}One can construct a signed permutation by inserting (barred or unbarred) elements
$n,n-1,\dots,1$ one at a time from largest to smallest.  To construct a signed permutation in which
the right-to-left maxima are unbarred, first we insert
$n$ into the empty word.  There is only one way to do this since $n$ must be unbarred. Now suppose we
have inserted
$n,n-1,\dots,n-j+1$.   If we insert
$n-j$ at the end of the partially constructed word, $n-j$ will be a right-to-left maximum in the final
word and must be unbarred.  If we insert
$n-j$ in any of the other $j$ positions $n-j$ will not be a right-to-left maximum and can therefore
be barred or unbarred.  Hence there are $2j+1$ ways to insert $n-j$, for each $j = 0,1,\dots,n-1$.
\end{proof}

\begin{Rem} \label{dow}{\rm The posets $\Pi_n$ and $\Pi^B_n$ are examples of Dowling
lattices. Gottlieb and Wachs
\cite[Section 9]{GW} have constructed splitting bases for general Dowling lattices. 
The splitting basis for $\Pi_n$ is a special case of this Dowling lattice construction.  For
$\Pi_n^B$, the  Dowling lattice construction produces the  basis consisting of 
cycles $\rho_{\omega,\epsilon}$ where
$(\omega,\epsilon) \in B_n$ is such that all {\it left-to-right} (rather than right-to-left) maxima
of
$(\omega,\epsilon)$ are unbarred.  Although this basis is similar in appearance to the
type $B$ splitting basis given here, one can show that it cannot be obtained from our geometric
construction.  (Consider   the intersection of $\mathcal B_3$ with the cube as in Figure 4.  The
union of the closure of the regions corresponding to the signed permutations whose left-to-right
maxima are unbarred, is not  simply connected. Hence these regions cannot be the bounded regions
of  a generic slice.)  However by using techniques similar to that of
\cite{GW}, one can get a variation of the Dowling lattice splitting basis which does reduce to the
type $B$ splitting basis  given here.  

Although the only Dowling lattices that are
intersection lattices of real hyperplane arrangements are the partition lattice and the signed
partition lattice, there  are other
 Dowling lattices that are intersection lattices of complex hyperplane
arrangements (cf.
\cite[Section 8]{GW}).  It would be interesting to find a geometric interpretation of the
above mentioned variation of the splitting basis for such Dowling lattices.  }
\end{Rem}

\begin{Rem}{\rm The splitting basis  of type $A$ is used in \cite{Wa96} and \cite{Wa98} to
obtain information about  the representation of the symmetric group on the homology of the
partition lattice.  Unfortunately, the
splitting basis of type
$B$ (or
$D$) does not appear to reveal much about the representation of the Coxeter group on  homology.
See \cite{Ha84}, \cite{Be91} and \cite{GW} for work on the representation of $B_n$ on the
homology of the signed partition lattice.  
 }
\end{Rem}

\section{Type $D$}

The type $D$ Coxeter arrangement is the  hyperplane arrangement $$ \DD_n =\{ x_i = 
x_j \mid 1
\le i < j
\le n\} \, \cup \, \{ x_i =  - x_j \mid 1 \le i < j
\le n\} . $$  The orthogonal reflections across
 the hyperplanes generate the Coxeter group $D_n$, which is the subgroup of $B_n$ consisting of
all signed permutations that have an even number of bars.  Clearly $\DD_n$ is a
subarrangement
of $\BB_n$ and its intersection lattice is isomorphic to  $\Pi_n^D$, the join-sublattice of
$\Pi^{B}_n$ consisting of all signed partitions  whose zero block does not have size
$2$. 
This isomorphism is denoted by $\gamma: L_{{\mathcal D}_n} \to \Pi_n^D $ and is the restriction of
the isomorphism
$\gamma: L_{\BB_n} \to \Pi_n^B $ defined in the previous section.   
  
The arrangement $\DD_n$ is essential and has $2^{n-1} n!$ regions which are simplicial cones in 
one-to-one correspondence with the elements of $D_n$.  Under this correspondence the signed
permutation
$(\omega,\epsilon)$ corresponds to the region
$$\wt R_{\omega,\epsilon} = \{\bold  x \in \R^n \mid \,\,  |x_{\omega(1)} | < \epsilon_2
x_{\omega(2)}<
\cdots < \epsilon_n x_{\omega(n)}\}.$$
Note that the hyperplane $x_{w(1)} = 0$ divides the region $\wt R_{\omega,\epsilon}$ into the
regions
$R_{\omega,\epsilon}$
and  $ R_{\omega,\epsilon^\prime}$,
where $\epsilon^\prime_i = \epsilon_i$ for $i = 2,3,\dots,n$ and $\epsilon^\prime_1 = - \epsilon_1$.

For each signed permutation $(\omega,\epsilon) \in D_n$, let $\wt\Pi_{\omega,\epsilon}$ be the
induced subposet of $\Pi^D_n$ consisting of all signed partitions obtained by splitting either
the signed permutation $(\omega,\epsilon)$ or the signed permutation $(\omega,\epsilon^\prime)$ at
all positions in a subset of
$\{0,1,\dots,n-1\}$ whose smallest element is not $1$.  
Although it is not as evident as for types $A$ and
$B$, this induced subposet is also a subset lattice.

\begin{Prop} The induced subposet
$
\wt
\Pi_{\omega,\epsilon}$ of $\Pi_n^D$ is isomorphic to the lattice of subsets of $[n]$.  
\end{Prop}

\begin{proof}
 Define the map $f: 2^{\{0,1,\dots,n-1\}} \to \wt
\Pi_{\omega,\epsilon}$ by letting $f(S)$ be the signed partition obtained by splitting 
$(\omega,\epsilon)$ at all
positions in $S$ if $1$ is not the smallest element of $S$, and by splitting 
$(\omega,\epsilon^\prime)$ at all positions in $S-\{1\} \cup \{0\}$  otherwise. We leave it to the
reader to check that this map is an isomorphism from the dual of $2^{\{0,1,\dots,n-1\}}$ to 
$\wt
\Pi_{\omega,\epsilon}$.  
\end{proof}

\begin{Prop}\label{fundd} For all $(\omega,\epsilon) \in D_n$, the image $\gamma(\rho_{\wt
R_{\omega,\epsilon}})$ is the fundamental cycle of $\overline{\wt \Pi}_{\omega,\epsilon}$. 
\end{Prop}

\begin{proof}  We show that $\gamma$ restricts to an isomorphism from the subposet $z(P_{\wt
R_{\omega,\epsilon}})$ of $L_{\DD_n}$ to $\wt \Pi_{\omega,\epsilon}$.  The elements of $P_{\wt
R_{\omega,\epsilon}}$ are sets that have one of  the following forms:
\begin{eqnarray*}\{\bold x\in \R^n \mid 0= x_{\omega(1)} = \dots = x_{\omega(i_1)} <
\epsilon_{i_1+1} x_{\omega(i_1+1)} =
\dots = \epsilon_{i_2} x_{\omega(i_2)}
 < \cdots \\
\cdots < \epsilon_{i_{k}+1} x_{\omega(i_{k}+1)} = \dots =  \epsilon_{n}
x_{\omega(n)}\},\end{eqnarray*}
\begin{eqnarray*}\{\bold x\in \R^n \mid0 < x_{\omega(1)} = \epsilon_2 x_{\omega(2)} = \dots =
\epsilon_{i_1}x_{\omega(i_1)} <
\epsilon_{i_1+1} x_{\omega(i_1+1)} =
\dots = \epsilon_{i_2} x_{\omega(i_2)}
\\ < 
\cdots < \epsilon_{i_{k}+1} x_{\omega(i_{k}+1)} = \dots =  \epsilon_{n}
x_{\omega(n)}\},\end{eqnarray*}
\begin{eqnarray*}  \{\bold x\in \R^n \mid 0 < -x_{\omega(1)} = \epsilon_2 x_{\omega(2)} = \dots =
\epsilon_{i_1}x_{\omega(i_1)} <
\epsilon_{i_1+1} x_{\omega(i_1+1)} =
\dots = \epsilon_{i_2} x_{\omega(i_2)}
\\ < 
\dots < \epsilon_{i_{k}+1} x_{\omega(i_{k}+1)} = \dots =  \epsilon_{n}
x_{\omega(n)}\},\end{eqnarray*}
\begin{eqnarray*}  \{\bold x\in \R^n \mid 0 <|x_{\omega(1)}| < \epsilon_2 x_{\omega(2)} = \dots =
\epsilon_{i_1}x_{\omega(i_1)} <
\epsilon_{i_1+1} x_{\omega(i_1+1)} =
\dots = \epsilon_{i_2} x_{\omega(i_2)}
\\ < 
\dots < \epsilon_{i_{k}+1} x_{\omega(i_{k}+1)} = \dots =  \epsilon_{n}
x_{\omega(n)}\},\end{eqnarray*} where
$2 \le i_1 < \dots < i_{k} < n$.

The linear span of the cells are the respective subspaces:
\begin{eqnarray*}\{\bold x\in \R^n \mid 0= x_{\omega(1)} = \dots = x_{\omega(i_1)} ,\,
\epsilon_{i_1+1} x_{\omega(i_1+1)} =
\dots = \epsilon_{i_2} x_{\omega(i_2)}
 , \, \dots \\
\dots ,\, \epsilon_{i_k+1} x_{\omega(i_k+1)} = \dots =  \epsilon_{n}
x_{\omega(n)}\},\end{eqnarray*}
\begin{eqnarray*}\{\bold x\in \R^n \mid  x_{\omega(1)} = \epsilon_2 x_{\omega(2)} = \dots =
\epsilon_{i_1}x_{\omega(i_1)} ,\,
\epsilon_{i_1+1} x_{\omega(i_1+1)} =
\dots = \epsilon_{i_2} x_{\omega(i_2)},\dots
\\  
\dots, \epsilon_{i_k+1} x_{\omega(i_k+1)} = \dots =  \epsilon_{n}
x_{\omega(n)}\},\end{eqnarray*}
\begin{eqnarray*}  \{\bold x\in \R^n \mid -x_{\omega(1)} = \epsilon_2 x_{\omega(2)} = \dots =
\epsilon_{i_1}x_{\omega(i_1)} ,\,
\epsilon_{i_1+1} x_{\omega(i_1+1)} =
\dots = \epsilon_{i_2} x_{\omega(i_2)}, \dots
\\ 
\dots , \epsilon_{i_k+1} x_{\omega(i_k+1)} = \dots =  \epsilon_{n}
x_{\omega(n)}\},\end{eqnarray*}
\begin{eqnarray*}  \{\bold x\in \R^n \mid \epsilon_2 x_{\omega(2)} = \dots =
\epsilon_{i_1}x_{\omega(i_1)} ,\,
\epsilon_{i_1+1} x_{\omega(i_1+1)} =
\dots = \epsilon_{i_2} x_{\omega(i_2)}, \dots
\\  
\dots, \epsilon_{i_k+1} x_{\omega(i_k+1)} = \dots =  \epsilon_{n}
x_{\omega(n)}\},\end{eqnarray*} where
$2 \le i_1 < \dots < i_k < n$.

Clearly the images of such subspaces under $\gamma$ are precisely the elements of $\wt
\Pi_{\omega,\epsilon}$.
\end{proof}

\begin{Prop}\label{boundd} Let $ \bold v  = (1,2,2^2,\dots,2^{n-1})$.  Then the affine hyperplane $H_\bold v$
is generic  with respect to the arrangement $\DD_n$.  Moreover, for all $(\omega,\epsilon) \in
D_n$,
$\wt R_{\omega,\epsilon}
\cap H_\bold v$ is bounded if and only if $\omega(1) \ne n$ and all right-to-left maxima of
$(\omega,\epsilon)$ are unbarred.
\end{Prop}

\begin{proof} Since by Proposition~\ref{bounded}, $H_\bold v$ is generic with respect
to the arrangement $\BB_n$,  it is generic with respect to any
subarrangement of $\BB_n$; in particular it is generic with respect to $\DD_n$.

   Since the hyperplane $x_{\omega(1)} = 0$ divides
$\wt R_{\omega,\epsilon}$ into the regions $ R_{\omega,\epsilon}$ and $
R_{\omega,\epsilon^\prime}$, the region $\wt R_{\omega,\epsilon}\cap H_\bold v$ is bounded if and
only if both
$ R_{\omega,\epsilon}\cap H_\bold v$ and $ R_{\omega,\epsilon^\prime}\cap H_\bold v$ are bounded.  By
Proposition~\ref{bounded}, both regions are bounded if and only if all right-to-left maxima of
both $(\omega,\epsilon)$ and  $(\omega,\epsilon^\prime)$ are unbarred.  This happens if and only if
$\omega(1)$ is not a right-to-left maximum, i.e., $\omega(1) \ne n$, and all right-to-left maxima of
$(\omega,\epsilon)$ are unbarred.  
\end{proof}

\begin{Thm}[type D splitting basis] \label{splitd}  For each $(\omega,\epsilon) \in D_n$, let
$\wt\rho_{\omega,\epsilon}$ be the fundamental cycle of $\overline{\wt \Pi}_{\omega,\epsilon}$. 
Then 
$$\{\wt\rho_{\omega,\epsilon} \mid 
(\omega,\epsilon) \in D_n,\,\, \omega(1) \ne n  \text { and all right-to-left maxima of}$$
$$(\omega,
\epsilon) \text{ are unbarred} 
\} $$
is a basis for  $\wt H_{n-2}(\overline{\Pi^D_n})$.
\end{Thm}

\begin{proof} 
This follows from Theorem \ref{T2} and Propositions~\ref{fundd} and~\ref{boundd}.
\end{proof}

The following corollary is well-known in the theory of Coxeter arrangements.
\begin{Cor}
  The rank of $\wt H_{n-2}(\overline{\Pi_n^D})$ is $$1 \cdot 3 \cdot 5 \cdots
(2n-3)
\cdot (n-1).$$
\end{Cor}

\begin{proof} We construct a signed permutation $(\omega,\epsilon) \in D_n$ with all right-to-left
maxima unbarred and $\omega(1) \ne n$ by first choosing $\omega(1)$ in $n-1$ ways and then choosing
the signed permutation $(\omega(2)\omega(3)\dots \omega(n), \epsilon_2\epsilon_3\cdots\epsilon_n)$ so
that every right-to-left maximum is unbarred.  By the proof of Corollary~\ref{Dow}, this signed
permutation on $n-1$ letters can be chosen in  $1\cdot 3 \cdots (2n-3)$ ways.  If the number
of bars in $(\omega(2)\omega(3)\dots \omega(n), \epsilon_2\epsilon_3\cdots\epsilon_n)$ is even then
$\omega(1)$ must be unbarred; otherwise $\omega(1)$ must be barred. 
\end{proof}

\section{Interpolating partition lattices}

We now consider a family of posets which interpolates between the type D partition lattice and the
type B partition lattice.  For $T \subseteq [n]$, let $\Pi^{DB}_n(T)$ be the join-sublattice of
$\Pi^{B}_n$ consisting of all signed  partitions whose zero block is not 
$\{0,a\}$  for
$a
\in [n]\setminus T$. Clearly
$\Pi^{DB}_n(\emptyset) = \Pi^{D}_n$ and $\Pi^{DB}_n([n]) = \Pi^B_n$.  The lattice $\Pi^{DB}_n(T)$ is
the intersection lattice of the hyperplane arrangement
 \begin{eqnarray*} \DD\BB_n(T)  & = &\{ x_i = 
x_j \mid 1
\le i < j
\le n\} \, \cup \, \{ x_i =  - x_j \mid 1 \le i < j
\le n\}   \\ & &\, \cup \,\,\, \{x_i = 0 \mid i \in T\}.\end{eqnarray*}
These interpolating arrangements were introduced by Zaslavsky \cite{Za81}.
The following theorem 
generalizes Theorems~\ref{splitb} and~\ref{splitd}.

\begin{Thm}\label{DB} For each $(\omega,\epsilon) \in B_n$,  let
$\rho_{\omega,\epsilon}$ be the fundamental cycle of $\overline{ \Pi}_{\omega,\epsilon}$, and
if $(\omega , \epsilon) \in D_n $ let
$\wt\rho_{\omega,\epsilon}$ be the fundamental cycle of $\overline{\wt \Pi}_{\omega,\epsilon}$.
The set  
\begin{eqnarray*}\{\rho_{\omega,\epsilon} \mid & &\!\!\!\! \!\!\!\!
(\omega,\epsilon) \in B_n,\,\,  \omega(1)  \in  T \text{ and  all right-to-left  maxima  of } \\
& &(\omega,
\epsilon) \text{ are unbarred} 
\} \,\,\, \\ & & \cup \\  \{\wt\rho_{\omega,\epsilon} \mid & & \!\!\!\!\!\!\!
(\omega,\epsilon) \in D_n,\,\, \omega(1) \notin T \cup \{n\}    \text{ and all right-to-left
maxima  of } \\ & & (\omega,
\epsilon) \text{ are unbarred} 
\} \end{eqnarray*}
forms a basis for $\wt H_{n-2}(\overline{\Pi^{DB}_n(T)})$.
\end{Thm}

\begin{proof}   There are two types of regions
of the hyperplane arrangement $\DD\BB_n(T)$, namely $R_{\omega,\epsilon}$  for $(\omega,\epsilon) \in 
B_n$ and $\omega(1) \in T$, and 
$\wt R_{\omega,\epsilon}$ for $(\omega,\epsilon) \in 
D_n$ and $\omega(1) \notin T$.  By Propositions~\ref{fundb} and~\ref{fundd}, 
$\rho_{\omega,\epsilon}$ and
$\wt \rho_{\omega,\epsilon}$ are the respective images (up to sign) of
$\rho_{R_{\omega,\epsilon}}$ and $\rho_{\wt R_{\omega,\epsilon}}$,  under the isomorphism
$\gamma$. 

We proceed as in the proof of Theorem~\ref{splitd}.  For $ \bold v  = (1,2,\dots,2^{n-1})$, 
$H_\bold v$ is generic with respect to the arrangement $\DD\BB_n(T)$.  Hence 
the result now follows from Theorem~\ref{T2}, Propositions~\ref{bounded} and~\ref{boundd}.
 \end{proof}

\begin{Cor}[Jambu and Terao {\cite{JT}}] The rank of $\wt
H_{n-2}(\overline{\Pi_n^{DB}(T)})$ is
$$1
\cdot 3
\cdot 5
\cdots (2n-3)
\cdot (|T|+n-1).$$
\end{Cor}

\begin{proof} Suppose $n\in T$.  Then the number of signed permutations $(\omega,\epsilon) \in B_n$
such that $\omega(1) \in T$ and all right-to-left maxima are unbarred is $(2|T|-1) \cdot 1 \cdot 3
\cdots (2n-3)$.  The number of signed permutations $(\omega,\epsilon) \in D_n$ such that $\omega(1)
\notin T \cup
 \{n\}$ and all right-to-left maxima are unbarred is $(n-|T|) \cdot1 \cdot 3 \cdots
(2n-3)$.  

Now suppose $n \notin T$. Then the number of signed permutations $(\omega,\epsilon) \in B_n$ such
that $\omega(1) \in T$ and all right-to-left maxima are unbarred is $2|T| \cdot 1 \cdot 3 \cdots
(2n-3)$.  The number of signed permutations $(\omega,\epsilon) \in D_n$ such that $\omega(1) \notin T
\cup \{n\}$ and all right-to-left maxima are unbarred is $(n-|T|-1) \cdot1 \cdot 3 \cdots
(2n-3)$.  In either case, the total number of elements in the basis is $(|T|+n-1) \cdot1 \cdot 3
\cdots (2n-3)$.
\end{proof}

J\'ozefiak and Sagan \cite{JS} have studied other families of hyperplane arrangements which
interpolate between Coxeter arrangements.  One can apply our results to these arrangements.  
The family of arrangements that interpolate between $\Aa_{n-2}$ and $\Aa_{n-1}$ are particularly
amenable to our approach.  For $T \subseteq [n-1]$, let $\Aa_n(T)$ be the arrangement in $\R^n$,
$$ \Aa_n(T)  = \{ x_i = 
x_j \mid 1
\le i < j
\le n-1\} \, \cup \, \{x_n = x_i \mid i \in T\}.$$

Let $\Pi_n(T)$ be the induced subposet of $\Pi_n$
consisting of all partitions $\pi$ such that the block of $\pi$ containing $n$ is either a singleton
or has nonempty intersection with $T$. This is the intersection lattice  of $\Aa_n(T)$.  The following
theorem generalizes Theorem~\ref{splita} by providing a splitting basis for the homology of
$\Pi_n(T)$.
  Recall that in Section~\ref{braid} we defined   $\rho_\omega$ to be the fundamental cycle of
$\overline{
\Pi}_\omega$, for each
$\omega\in S_n$.

\begin{Thm}  \label{AT} For $\emptyset \ne T \subseteq [n-1]$,
the set  
$$\{\rho_\omega \mid 
\omega \in S_n,\,\,  \omega(n) =n\text{ and }  \omega(n-1) \in T\}$$
forms a basis for $\wt H_{n-3}(\overline{\Pi_n(T)})$.
\end{Thm}

\begin{proof} The proof, which is similar to  that of Theorem~\ref{DB}, uses results from
Section~\ref{braid} and is left to the reader.
\end{proof}

\begin{Rem} {\rm
There are easier and more direct ways to prove Theorem~\ref{AT}.  For instance, one can restrict the
second EL-labeling for
$\Pi_n$ given in the proof of Theorem~6.2 of \cite{Wa96}, to $\Pi_n(T)$. The  induced
shelling basis is precisely the basis given in Theorem~\ref{AT}.  }
\end{Rem}

\begin{Cor}[J\'ozefiak and Sagan \cite{JS}] For $\emptyset \ne T \subseteq [n-1]$, the rank of $\wt
H_{n-3}(\overline{\Pi_n(T)})$ is
$(n-2)!
\cdot |T|$.
\end{Cor}

In  \cite{St82}, Stanley showed that the restriction to $S_{n-1}$ of the representation of $S_n$
on
$\wt H_{n-3}(\Pi_n,\C)$ is the regular representation. In \cite{Wa96}, it was observed that the
splitting basis for
$\wt H_{n-3}(\overline{\Pi_n})$ makes this fact
transparent.  Indeed, the permutations that fix $n$ permute the basis cycles $\rho_\omega$. A
similar phenomenon occurs for $\Pi_n(T)$.

\begin{Cor}\label{reg} For $\emptyset \ne T\subseteq [n-1]$,  the representation
of 
$S_T \times S_{[n-1]\setminus T}$ on $\wt
H_{n-3}(\overline{\Pi_n(T)},\C)$ is isomorphic to the direct sum of $\binom {n-2} {|T| -1}$ copies of
the regular representation of $S_T \times S_{[n-1]\setminus T}$.
\end{Cor}

\begin{proof} If $\omega \in S_n$ satisfies $\omega(n) = n$ and $\omega(n-1) \in T$ then so does
$\sigma \omega$ for all $\sigma \in S_T \times S_{[n-1]\setminus T} \times S_{\{n\}}$.   In fact,
the elements of $T$ occupy the same set of positions in $\sigma\omega$ as in $ \omega$.  
Since
$\sigma
\rho_\omega = \rho_{\sigma\omega}$, we see that  $S_T \times S_{[n-1]\setminus T}
$ acts on $\wt
H_{n-3}(\overline{\Pi_n(T)},\C)$ by permuting basis cycles.  Also, the  orbit of $\rho_\omega$
is determined by the set of positions  that elements of $T$ occupy in $\omega$.  Hence the number of
orbits is $\binom {n-2} {|T| -1}$.  
\end{proof}

\begin{Rem}{\rm John Shareshian [personal communication] has found an alternative  proof of
Corollary~\ref{reg} which involves a computation of the M\"obius function of the
$\sigma$-invariant subposet of
$\Pi_n(T)$ for 
$\sigma \in S_T \times S_{[n-1]\setminus T}$ by means of Crapo's complementation formula
\cite{Cr66}.  }  
\end{Rem}

\begin{Rem} {\rm Another class of partition posets with a splitting basis is  the
class of $d$-divisible partition lattices \cite{Wa96}, or the more general restricted block size
partition lattices considered in \cite{Br} and \cite{BrW}.  These are not geometric lattices in
general; but they are
 intersection lattices of subspace arrangements and they are shellable.  The cycles in
the basis are polytopal.   The symmetric group acts on these lattices and the splitting basis
reveals much information about the representation of the symmetric group on homology. It
would be interesting to find a geometric explanation for this spitting basis.   One might also
consider    restricted block size partition subposets of $\Pi_n(T)$.  
 }\end{Rem}

\end{document}